\documentclass{tran-l}
\input{psfig}
\usepackage{amssymb}
 
 \newcommand{\C}{{\mathbb C}}
 
\newcommand{\X}{{\mathfrak X}}
\numberwithin{equation}{section}
 \newtheorem{theorem}{Theorem}[section]
 \newtheorem{lemma}[theorem]{Lemma}
 \newtheorem{proposition}[theorem]{Proposition}
 \newtheorem{corollary}[theorem]{Corollary}
 
 \newtheorem{example}[theorem]{Example}
 \newtheorem{definition}[theorem]{Definition}

 \newtheorem{fact}[theorem]{Fact}
\newcommand{\cal}{\mathcal}
 \def\boxx
  {\hfill \thinspace\vbox{\hrule height .5pt \hbox{\vrule
   width .5pt \vbox to 7pt{\hbox to 3.5pt{}} \vrule width .5pt}
   \hrule height 0pt depth .5pt}}

\begin{document}

\title{$SL_n$-Character Varieties as Spaces of Graphs}

\author{Adam S. Sikora}

\address{Department of Mathematics, University of Maryland,
College Park, MD 20742}

\email{asikora@math.umd.edu}
\thanks{Partially supported by the NSF grant DMS93-22675.}


\date{May 12, 1999 and, in revised form, June 5, 2000.}


\begin{abstract} An $SL_n$-character of a group $G$ is the trace of an
$SL_n$-representation of $G.$ We show that all algebraic relations
between $SL_n$-characters of $G$ can be visualized as relations between graphs 
(resembling Feynman diagrams) in any topological space $X,$ with
$\pi_1(X)=G.$ We also show that all such relations are implied by a
single local relation between graphs. In this way, we provide a
topological approach to the study of $SL_n$-representations of groups.

The motivation for this paper was our work with J. Przytycki
on invariants of links in 3-manifolds which are based on the
Kauffman bracket skein relation. These invariants lead to a notion of a
skein module of $M$ which, by a theorem of Bullock, Przytycki, and the
author, is a deformation of the $SL_2$-character variety of $\pi_1(M).$
This paper provides a generalization of this result to all $SL_n$-character
varieties.
\end{abstract}

\maketitle


\section{Introduction}
In this paper we introduce a new method in the study of representations of
groups into affine algebraic groups.  Although we consider only
$SL_n$-representations, the results of this paper can be generalized
to other affine algebraic groups, see \cite{Si-char}.

For any group $G$ and any commutative ring $R$ with $1$ there
is a commutative $R$-algebra $Rep_n^R(G)$ and {\it the universal
$SL_n$-representation}
$$j_{G,n}:G\to SL_n(Rep_n^R(G))$$ such that any representation
of $G$ into $SL_n(A),$ where $A$ is an $R$-algebra, factors through
$j_{G,n}$ in a unique way. This universal property uniquely 
determines $Rep_n^R(G)$ and $j_{G,n}$ up to an
isomorphism. 

$GL_n(R)$ acts on $Rep_n^R(G)$ (see Section 2) and the subring of
$Rep_n^R(G)$ composed of the elements fixed by the action,
$Rep_n^R(G)^{GL_n(R)},$ is called the {\it the universal
$SL_n$-character ring of $G$.}
This ring contains essential information about
$SL_n$-representations of $G.$ In particular, if $R$ is an
algebraically closed field of characteristic $0$ then there are natural
bijections between the following three sets:
\begin{itemize}
\item the set of all $R$-algebra homomorphisms
$Rep_n^R(G)^{GL_n(R)}\to R$
\item the set of all semisimple $SL_n(R)$-representations of $G$ up to
conjugation
\item the set of $SL_n(R)$-characters of $G.$
\end{itemize}

It is convenient to think about $Rep_n^R(G)^{GL_n(R)}$ as the
coordinate ring of a scheme, $\X_n(G)=Spec(Rep_n^R(G)^{GL_n(R)}),$
called the $SL_n$-character variety of $G.$
As explained in Section 6, the algebra $Rep_n^R(G)^{GL_n(R)}$ encodes
all algebraic relations between the $SL_n$-characters of $G.$
Unfortunately, it is very difficult to give a finite
presentation of $Rep_n^R(G)^{GL_n(R)}$ and, hence, to describe
$\X_n(G),$ even for groups $G$ with
relatively simple presentations.

In this paper, we present a topological approach to the study
of $SL_n$-character varieties. We prove that
$R[\X_n(G)]=Rep_n^R(G)^{GL_n(R)}$ is spanned by a special class of graphs
(resembling Feynman diagrams)
in $X,$ where $X$ is any topological space with $\pi_1(X)=G,$
see Theorem \ref{2.7}. Moreover, all relations between
the elements of this spanning set are induced by specific local
relations between the graphs, called {\it skein relations}.

We postpone a detailed study of applications of our graphical calculus
to the theory of $SL_n$-representations of groups to future papers.
In this paper, we content ourself with an example, in which we apply
our method to a study of $SL_3$-representations of the free group on
two generators. In this algebraically nontrivial example a huge
reduction of computational difficulties can be achieved by the
application of our geometric method.

This work is related to several areas of mathematics and physics:

{\bf Knot theory (and skein modules)} Skein relations between
links were used to define the famous polynomial invariants of
links, like the Conway, Jones, and Homfly polynomials, \cite{Co,
Jo, FYHLMO, P-T, Ka}.
In this paper we apply skein relations to the representations of
groups.

The motivation for this work was our earlier work on skein
modules, \cite{PS-2}. The main theorem of this paper generalizes
the Bullock-Przytycki-Sikora theorem relating the Kauffman bracket
skein module of a manifold $M$ to the $SL_2(\C)$-character variety
of $\pi_1(M),$ see \cite{B-2,PS-2}.

{\bf Quantum invariants of $3$-manifolds}
We hope that this work will help understand the 
connections between quantum invariants of $3$-manifolds
and representations of their fundamental groups.
It follows from the work of Yokota \cite{Yo} that for any
$3$-manifold $M,$ the $SU_n$-quantum
invariants of $M$ can be defined by using our graphs considered up to
relations which are $q$-deformations of our skein relations.

{\bf Spin networks and Gauge Theory} The graphs considered in this paper
have an interpretation as spin networks, see \cite{Si-char}.
They are also very similar to graphs used by physicists in non-abelian 
gauge theory (QCD),
see \cite{Cv}.

{\bf Number theory} After a preliminary version of this paper was
made available, M. Kapranov pointed out to us, that our work is
related to the work Wiles and others on
``pseudo-representations.'' In his work (related to Fremat's Last 
Theorem), Wiles gave necessary and
sufficient conditions under which a complex valued function on $G$
is a $GL_2(\C)$-character of $G.$ His ideas were developed further
and generalized to all $GL_n$-characters by Taylor, \cite{Ta}. See
also \cite{Ny, Ro}. These results provide a description of 
the coordinate ring of $GL_n$-character varieties quotiented by
nilpotent elements. Our
results are similar in spirit, but they are concerned with
$SL_n$-representations and they are stronger, since they describe
$R[\X_n(G)]$ (ie. $Rep_n^R(G)^{GL_n(R)}$) exactly  
(with possible nilpotent elements).

The plan of this paper is as follows. In Section 2 we introduce
some basic notions and facts concerning representations of groups.
In Section 3 we define the algebra
${\mathbb A}_n(X)$ in terms of graphs in $X$ and formulate
(Theorems \ref{2.6} and \ref{2.7}) the main results of the paper
asserting that ${\mathbb A}_n(X)$ is isomorphic to
$Rep_n^R(G)^{GL_n(R)},$ where $G= \pi_1(X).$ The proof requires
introducing another algebra, ${\mathbb A}_n(X,x_0),$ associated with any pointed
topological space $(X,x_0).$ The algebra ${\mathbb A}_n(X,x_0)$ is
an interesting object by itself, and for $n=2$ it already appeared in
the theory of skein modules as a relative skein algebra.
Sections 4 and 5 are devoted to the proof of the results of Section 3.
In the final section we consider trace identities and use
our results to describe the $SL_3$-character variety of the free group
on two generators.

I would  like to thank C. Frohman, B. Goldman, and J. H. Przytycki
for many helpful conversations. C. Frohman explained to me some
ideas which I subsequently used in this paper. 
Additionally, I would like thank C. Frohman and W. Goldman for helping me 
to make this paper more ``reader-friendly.''


\section{Background from representation theory}

In this section we introduce the basic elements of the theory of
$SL_n$-representations of groups. We follow the approach of Brumfiel
and Hilden, (\cite{B-H}, Chapter 8), which although formally restricted to
$SL_2$-representations, has a straightforward generalization to
$SL_n$-representations for any $n.$ Compare also \cite{L-M}, \cite{Pro-2}.

For any group $G$ and a commutative ring $R$ with $1,$ there
is a commutative $R$-algebra $Rep_n^R(G),$ 
and {\it the universal $SL_n$-representation} $$j_{G,n}:G\to
SL_n(Rep_n^R(G)),$$ with the following property:

\noindent For any commutative $R$-algebra $A$ and any representation
$\rho: G\to SL_n(A)$ there is a unique homomorphism of
$R$-algebras $h_{\rho}:Rep_n^R(G)\to A$ which induces a homomorphism
of groups
$$SL_n(h_{\rho}):SL_n(Rep_n^R(G))\to SL_n(A)$$
such that the following diagram commutes:

\begin{center}
\begin{tabular}{ccc}
$G$ & $\stackrel{j_{G,n}}{\longrightarrow}$ & $SL_n(Rep_n^R(G))$\\
& $\stackrel{\rho}{\searrow}$ & $\Big\downarrow
\scriptstyle{SL_n(h_{\rho})}$\\
& & \hspace*{-0.5in} $SL_n(A)$\\
\end{tabular}
\end{center}

\noindent This universal property determines $Rep_n^R(G)$
uniquely up to an isomorphism of $R$-algebras.
The algebra $Rep_n^R(G)$
may also be constructed explicitly in the following
way. Let $<g_i,\ i\in I|r_j,\ j\in J>$ be a presentation of $G$
such that all relations, $r_j,$ are monomials in non-negative
powers of generators, $g_i.$ Such a presentation exists for every
group $G.$ Since we work with groups which are not necessarily
finitely presented, $I$ and $J$ may be infinite. Let $P_n(I)$ be
the ring of polynomials over $R$ in variables $x^i_{jk},$ where
$i\in I$ and $j,k\in \{1,2,...,n\}.$ Let $A_i,$ for $i\in I,$ be
the matrix $(x^i_{jk})\in M_n(P_n(I)).$ For any word
$r_j=g_{i_1}^{n_1}g_{i_2}^{n_2}...g_{i_k}^{n_k}$ consider the
corresponding matrix
$M_j=A_{i_1}^{n_1}A_{i_2}^{n_2}...A_{i_k}^{n_k} \in M_n(P_n(I)).$
Let ${\cal I}$ be the two-sided ideal in $P_n(I)$ generated by
polynomials $Det(A_i)-1,$ for $i\in I,$ and by all entries of
matrices $M_j-Id,$ for $j\in J,$ where $Id$ is the identity matrix.
We denote the quotient,
$P_n(I)/{\cal I},$ by $Rep_n^R(G)$ and the quotient map $P_n(I)\to
Rep_n^R(G)$ by $\eta.$ Let $\bar x^i_{jk}=\eta(x^i_{jk})$ and
$\bar A_i=(\bar x^i_{jk})\in M_n(Rep_n^R(G)).$

Note that we divided $P_n(I)$ by
all relations necessary for the existence of a representation
$$j_{G,n}: G\to SL_n(Rep_n^R(G))$$ such that $j_{G,n}(g_i)=\bar A_i.$
It is easy to see that the two above definitions of $Rep_n^R(G)$
are equivalent and that $j_{G,n}: G\to SL_n(Rep_n^R(G))$ is the
universal $SL_n$-representation of $G.$ 

Let $A\in GL_n(R).$ By the definition of 
$Rep_n^R(G),$ there is a unique homomorphism $f_A: Rep_n^R(G)\to
Rep_n^R(G),$ such that the following diagram commutes

\begin{center}
\begin{tabular}{ccc}
$G$ & $\stackrel{j_{G,n}}{\longrightarrow}$ & $SL_n(Rep_n^R(G))$\\
& $\stackrel{A^{-1}j_{G,n}A}{\searrow}$ & $\Big \downarrow
\scriptstyle{SL_n(f_A)}$\\
& & $SL_n(Rep_n^R(G))$\\
\end{tabular}
\end{center}

\noindent One easily observes that $f_A$ is an automorphism of
$Rep_n^R(G)$ and that the
assignment $A\to f_A$ defines a left action of $GL_n(R)$ on $Rep_n^R(G).$
We denote this action by $A*,$ i.e. $f_A(r)=A*r,$ for any $r\in
Rep_n^R(G)$ and $A\in GL_n(R).$
We call the ring $Rep_n^R(G)^{GL_n(R)}$ consisting of elements of
$Rep_n^R(G)$
fixed by the action of $GL_n(R)$ {\it the universal $SL_n$-character ring of
$G.$} This term indicates a connection between the ring $Rep_n^R(G)^{GL_n(R)}$
and $SL_n$-characters of $G,$ i.e. traces of $SL_n$-representations of
$G.$ In the simplest case, when $R$ is an algebraically closed field
of characteristic $0,$ and $G$ is finitely generated this connection
can be described as follows.
Every representation $\rho: G\to SL_n(R)$ induces a homomorphism
$h_{\rho}: Rep_n^R(G)\to R$ whose restriction to $Rep_n^R(G)^{GL_n(R)}$ we
denote by $h_{\rho}'.$  The next proposition follows from the
geometric invariant theory and from the results of \cite{L-M}.

\begin{proposition}\label{1.0}\ \\
Under the above assumptions, the following sets are in a natural
correspondence given by bijections $\rho \to h_{\rho}',$ and $\rho\to
\chi=tr\circ\rho:$
\begin{itemize}
\item The set of all semisimple $SL_n(R)$-representations of $G;$
\item The set of all $R$-homomorphisms $Rep_n^R(G)^{GL_n(R)}\to R;$
\item The set of all $SL_n(R)$-characters of $G.$
\end{itemize}
\end{proposition}

By the proposition, we can identify the above sets and denote them
by $X_n(G).$ By the second definition, $X_n(G)$ is the affine
algebraic set composed of the closed points of the
$SL_n$-character variety, $\X_n(G),$ defined in the introduction.
In other words,
$$R[X_n(G)]=R[\X_n(G)]/\sqrt{0}=Rep_n^R(G)^{GL_n(R)}/\sqrt{0}.$$
As shown in \cite{L-M, KM}, $Rep_n^R(G)^{GL_n(R)}$ may contain
nilpotent elements and, therefore $\X_n(G)$ contains more subtle
information about $SL_n$-representations of $G$ then $X_n(G).$ The
ring $Rep_n^R(G)^{GL_n(R)}$ will be given a topological
description in Section 3.

The $GL_n(R)$ action on $Rep_n^R(G)$ induces an action of $GL_n(R)$ on
the ring
of $n\times n$ matrices over $Rep_n^R(G).$ If $M=(m_{ij})\in M_n(Rep_n^R(G))$
and $A\in GL_n(R)$ then
\begin{equation}
A*M=A\left(
\begin{matrix}
A*m_{11} & A*m_{12} & ... & A*m_{1n}\cr \vdots &
\vdots & ... & \vdots\cr A*m_{n1} & A*m_{n2} & ... &
A*m_{nn}
\end{matrix}
\right)A^{-1}.\label{e1}
\end{equation}


There is an equivalent definition of the action of $GL_n(R)$ on $Rep_n^R(G)$ and
on $M_n(Rep_n^R(G)).$ In order to introduce it we will first define
$GL_n(R)$-actions on $P_n(I)$ and $M_n(P_n(I)).$
We can consider $P_n(I)$ as a ring of polynomial functions
defined on the product of $I$ copies of $M_n(R),\ M_n(R)^I\to R,$ by
identifying $x^{i_0}_{jk}\in P_n(I)$ with a map assigning to
$(M_i)_{i\in I}\in M_n(R)^I$ the $(j,k)$-entry of $M_{i_0}.$
Therefore, $$P_n(I)=Poly(M_n(R)^I,R).$$
With this identification any entry in a matrix $M$ in $M_n(P_n(I))$ is
a polynomial function on $M_n(R)^I.$ Therefore we can think of elements of
$M_n(P_n(I))$ as coordinate-wise polynomial functions
$M_n(R)^I\to M_n(R),$
$$M_n(P_n(I))=Poly(M_n(R)^I,M_n(R)).$$

If $X,Y$ are sets with a left
$G$-action, then the set of all functions $Fun(X,Y)$ has a natural
left $G$-action defined for any $f:X\to Y$ and $g\in G$ by
$g*f(x)=gf(g^{-1}x),$ for $x\in X.$
$GL_n$ acts on $M_n(R)$ and on $M_n(R)^I$ by conjugation and it acts
trivially on $R.$ These actions induce $GL_n(R)$-actions on
$Fun(M_n(R)^I,R)$ and on $Fun(M_n(R)^I,M_n(R)),$ which restrict to
$P_n(I)=Poly(M_n(R)^I,R)$ and
$M_n(P_n(I))=Poly(M_n(R)^I,M_n(R)).$
The following statement is a consequence of the above definitions.

\begin{lemma}\label{1.1}\ \\
\vspace*{-0.2in}
\begin{enumerate}
\item The natural embedding of $P_n(I)$ into $M_n(P_n(I))$ as scalar
matrices is $GL_n(R)$-equivariant.
\item $A_{i_0}=(x^{i_0}_{jk})\in M_n(P_n(I))$ is
invariant under the action of $GL_n(R),$ for any $i_0\in I.$
\end{enumerate}
\end{lemma}

%

Now, we are going to show that $\eta: P_n(I)\to Rep_n^R(G),$ is
$GL_n(R)$-equivariant and hence the action of $GL_n(R)$ on $P_n(I)$
induces a $GL_n(R)$-action on $Rep_n^R(G)$ which coincides with the
$GL_n(R)$-action on $Rep_n^R(G)$ defined previously.

\begin{proposition}\label{1.2}\ \\
The following diagram commutes:

\begin{equation}
\begin{tabular}{ccc}
$M_n(P_n(I))$ & $\stackrel{M_n(\eta)}{\longrightarrow}$ & $M_n(Rep_n^R(G))$\\
$\Big \downarrow \scriptstyle{Tr}$ & & $\Big \downarrow \scriptstyle{Tr}$\\
\hspace*{-.1in} $P_n(I)$ &$\stackrel{\eta}{\longrightarrow}$ & $Rep_n^R(G)$\\
\end{tabular} \label{e2}
\end{equation}
and all maps appearing in it intertwine with the $GL_n(R)$-action.
\end{proposition}

\begin{proof}
Since the commutativity of the above diagram is obvious we will prove only
that the trace functions and homomorphisms $\eta,\ M_n(\eta)$ are
$GL_n(R)$-equivariant.

\begin{itemize}
\item The trace map $Tr:M_n(R)\to R$ is $GL_n(R)$-equivariant.
Therefore the induced map
$$Tr: M_n(P_n(I))=Poly(M_n(R)^I,M_n(R))\to Poly(M_n(R)^I,R)=P_n(I)$$
is also $GL_n(R)$-equivariant.

\item If $M=(m_{ij})\in M_n(Rep_n^R(G))$ and $A\in GL_n(R),$ then
$A*M$ is given by matrix (\ref{e1}) whose trace is
$Tr(A*M)=\sum_{i=1}^n A*m_{ii}=A*Tr(M).$ Therefore
$$Tr:M_n(Rep_n^R(G))\to Rep_n^R(G)$$ is $GL_n(R)$-equivariant.

\item Recall that $P_n(I)$ is generated by elements $x^i_{jk}.$
Therefore in order to prove that $\eta$ is $GL_n(R)$-equivariant
it is enough to show that $\eta(A*x^i_{jk})=A*\bar x^i_{jk},$
for any $i\in I,\ j,k\in\{1,2,...,n\},$ where $\bar
x^i_{jk}=\eta(x^i_{jk})\in Rep_n^R(G).$

For any $i_0\in I,\ j_{G,n}(g_{i_0})=\bar A_{i_0}\in SL_n(Rep_n^R(G)).$
By the definition of the $GL_n(R)$-action on $Rep_n^R(G),$
$A^{-1}j_{G,n}(g_{i_0})A$ is the matrix obtained from $j_{G,n}(g_{i_0})$ by
acting on all its entries by $A.$ Therefore
\begin{equation}
A*\bar x^{i_0}_{jk}= (j,k)-{\rm entry\ of\ } A^{-1}\bar A_{i_0}A.
\label{e3}
\end{equation}

Having described $A*\bar x^{i_0}_{jk},$ we need to give an explicit
description of $A*x^{i_0}_{jk}\in
P_n(I).$ Recall that we identified $x^{i_0}_{jk}$ with the map
$M_n(R)^I\to R$ assigning to $\{M_i\}_{i\in I}$ the
$(j,k)$-entry of $M_{i_0}.$ The definition of the
$GL_n(R)$-action on maps between $GL_n(R)$-sets implies that
$$(A*x^{i_0}_{jk})(\{M_i\}_{i\in I})=A*\left(x^{i_0}_{jk}
(A^{-1}*\{M_i\}_{i\in I})\right).$$
Since $GL_n(R)$ acts by simultaneous conjugation on $M_n(R)^I$ and it
acts trivially on $R,$ the right side of the above equation
is equal
to the $(j,k)$-entry of $A^{-1}M_{i_0}A.$ But the entries of $M_{i_0}$
are given by the values of functions $x^{i_0}_{jk}$ evaluated on
$\{M_i\}_{i\in I}.$ Therefore
\begin{equation}
A*x^{i_0}_{jk}= (j,k)-{\rm entry\ of\ } A^{-1}(x^{i_0}_{jk})A.
\label{e4}
\end{equation}

Finally, (\ref{e3}) and (\ref{e4}) imply that
$$\eta(A*x^{i_0}_{jk})= \eta\left({\rm the\ }(j,k){\rm -entry\ of\ }
A^{-1}A_{i_0}A\right)=$$
$${\rm the\ }(j,k){\rm -entry\ of\ }A^{-1}\bar A_{i_0}A=A*\bar x^{i_0}_{jk}.$$

\item We prove that $M_n(\eta)$ is equivariant.
Let $M=(m_{jk})\in M_n(P_n(I)).$
Notice that the definition of $GL_n(R)$-action on $M_n(P_n(I))$ implies
that $A*M=A\left(A*m_{jk}\right)A^{-1}.$
Therefore
$$M_n(\eta)(A*M)= M_n(\eta)(A(A*m_{jk})A^{-1})=
A(\eta(A*m_{jk}))A^{-1}.$$
Since $\eta$ is $GL_n(R)$-equivariant, the matrix on the right side of
the above equation is
$A(A*\eta(m_{jk}))A^{-1}=A*(\eta(m_{jk})).$
Therefore $M_n(\eta)$ is also $GL_n(R)$-equivariant.
\end{itemize}
\end{proof}

The above proposition implies that there exists a function
$$Tr: M_n(Rep_n^R(G))^{GL_n(R)}\to Rep_n^R(G)^{GL_n(R)}.$$
This function will be given a simple topological interpretation in the next
section.

\begin{proposition}\label{1.3}\ \\
The image of the universal $SL_n$-representation $j_{G,n}: 
G\to M_n(Rep_n^R(G))$
is invariant under the action of $GL_n(R).$
\end{proposition}

\begin{proof}
Since the elements $g_i$ generate $G,$
it is sufficient to show that $j_{G,n}(g_i)\in M_n(Rep_n^R(G))
^{GL_n(R)}.$
By Lemma \ref{1.1}(2) $A_i\in M_n(P_n(I))^{GL_n(R)}.$
The map $M_n(\eta)$ is equivariant. Therefore it takes the invariant
$A_i$ to the invariant $M_n(\eta)(A_i)=j_{G,n}(g_i).$
\end{proof}


\section{Skein Algebras}

In this section we assign to each path connected topological space $X$
a commutative
$R$-algebra ${\mathbb A}_n(X)$ and to each pointed path connected
topological space $(X,x_0)$ an $R$-algebra ${\mathbb A}_n(X,x_0).$ These
algebras encode
the most important information about the $SL_n$-representations of
$\pi_1(X,x_0).$ We will show that if $R$ is a field of characteristic $0$
(but not necessarily algebraically closed) then ${\mathbb A}_n(X)$ is
isomorphic to the universal $SL_n$-character ring, $Rep_n^R(G)^{GL_n(R)},$ where
$G=\pi_1(X,x_0),$ and ${\mathbb A}_n(X,x_0)$ is isomorphic to
$M_n(Rep_n^R(G))^{GL_n(R)}.$

We start with a definition of a graph which is the most suitable for
our purposes.
A {\it graph} $D=({\cal V},{\cal E},{\cal L})$ consists of a vertex-set
${\cal V},$ a set of oriented edges  ${\cal E},$ and a set of oriented loops
${\cal L}.$ Each edge $E\in {\cal E}$ has a beginning $b(E)\in {\cal
V}$ and an end $e(E)\in {\cal V}.$ Loops have neither beginnings nor ends.
If $b(E)=v$ or $e(E)=v$ then $E$ is {\it incident} with $v.$
For any $v\in {\cal V}$ all edges incident to $v$ are are ordered by
consecutive integers $1,2,...$ 
Therefore the beginning and the end of each edge is
assigned a number.

The sets ${\cal V},{\cal E},{\cal L}$ are finite.
We topologize each graph as a CW-complex.
The topology of a graph coincides with the topology of its edges
$E\simeq [0,1], E\in {\cal E},$ and its loops $L\simeq S^1, L\in {\cal L}.$
There is a natural notion of isomorphism of graphs.

Let ${\cal G}$ be a set of representatives of all isomorphism classes
of graphs defined above.
We say that a vertex $v$ is an $n$-valent source of a
graph $D$ if $n$ distinct edges of $D$ begins at $v$ and no edge ends at $v.$
Similarly, we say that $v$ is an $n$-valent sink of $D$ if $n$ distinct edges
end at $v$ and no edge begins at $v.$
Let ${\cal G}_n$ denote the set of all graphs in ${\cal G},$
all of whose vertices are either $n$-valent sources or $n$-valent sinks.
We assume that the empty graph $\emptyset$ is also an element of ${\cal G}_n.$
We denote the single loop in ${\cal G}_n,$ i.e. the connected graph
without any vertices, by $S^1.$
Let ${\cal G}_n'$ denote the set of all graphs $D\in {\cal G}$ such that
$D$ has one $1$-valent source, one $1$-valent sink and all other vertices
of $D$ are $n$-valent sources or $n$-valent sinks. We denote the
single edge in ${\cal G}_n',$ i.e. the connected graph without any
$n$-valent vertices by $[0,1].$

Let $X$ be a path connected topological space.
We will call any continuous map $f: D\to X,$ where $D\in {\cal G}_n,$
{\it a graph in $X.$} We identify two maps $f_1,f_2:D\to X$ if they are
homotopic. Let us denote the set of all graphs in $X$ by ${\cal G}_n(X).$
Similarly, we define ${\cal G}_n(X,x_0)$ to be the set of all maps
$f: D\to X\times [0,1],$ where $D\in {\cal G}_n'$
and $f$ maps the $1$-valent sink of $D$ to $(x_0,1)$ and the $1$-valent
source of $D$ to $(x_0,0).$ We identify maps which are homotopic relative
to $(x_0,0)$ and $(x_0,1).$ We will call elements of
${\cal G}_n(X,x_0)$ {\it relative graphs in $X\times [0,1].$}

We introduce a few classes of graphs in ${\cal G}_n(X)$ and ${\cal G}_n(X,x_0)$
which will be often used later on in the paper.
Let $L_\gamma: S^1\to X$ be a graph in $X$ which represents the
conjugacy class of $\gamma \in \pi_1(X,x_0).$
We denote by $E_\gamma$ a relative graph $E_\gamma:[0,1]\to X\times [0,1],
\ E_\gamma(0)=(x_0,0),\ E_\gamma(1)=(x_0,1),$ whose projection into $X,$
$$[0,1]\stackrel{E_{\gamma}}{\to} X\times [0,1]\to X$$
represents $\gamma \in \pi_1(X,x_0).$ Let $EL_\gamma:[0,1]\cup S^1\to X\times
[0,1]$ be a relative graph such that
$EL_{\gamma}(t)=(x_0,t),$ for $t\in [0,1],$ and $EL_{\gamma|S^1}:
S^1\to X\times [0,1]\to X$ represents the conjugacy class of $\gamma\in
\pi_1(X,x_0).$

For any two graphs $f_1:D_1\to X$ and $f_2:D_2\to X,\ f_1,f_2\in
{\cal G}_n(X),$ we define a product of them to be
$f_1\cup f_2:D_1\cup D_2\to X,\ f_1\cup f_2 \in {\cal G}_n(X),$ where
$D_1\cup D_2$ denotes the disjoint union of $D_1$ and $D_2.$
Therefore the free $R$-module  $R{\cal G}_n(X)$ on ${\cal G}_n(X)$
can be considered as a commutative
$R$-algebra. The empty graph $\emptyset:\emptyset\to X$ is an identity in
$R{\cal G}_n(X).$

In the next definition we will represent fragments of diagrams by
coupons, as depicted below.

\centerline{\psfig{figure=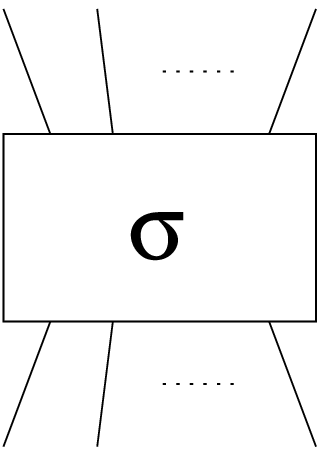,height=.7in}}

\noindent This coupon means a braid corresponding to a permutation
$\sigma\in S_n.$
\footnote{ Since we consider graphs up to homotopy equivalence,
it does not matter which braid corresponding to $\sigma$ we take.}

\begin{example}\label{2.1}\ \\
If $\sigma=(1,2,3)\in S_3$ then\\
\centerline{\psfig{figure=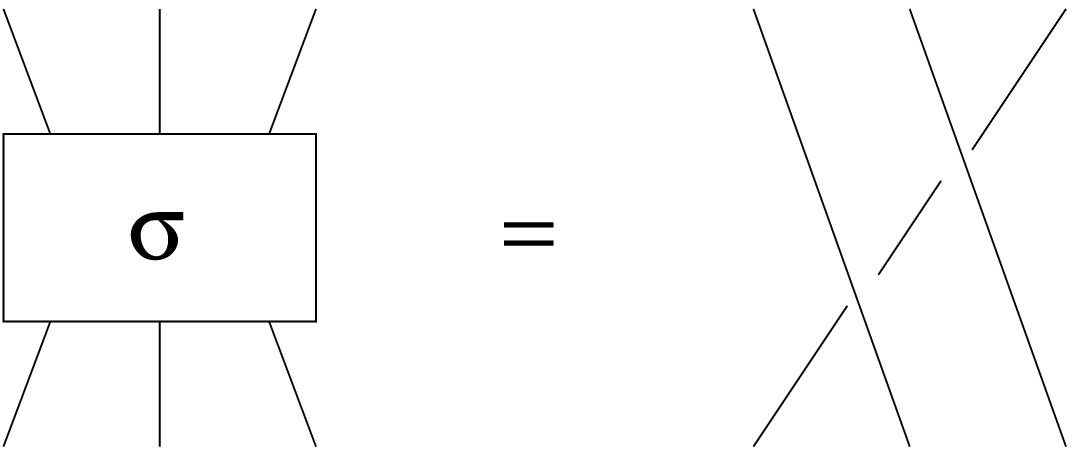,height=.7in}}
\end{example}

Suppose that $f:D\to X$ is a graph in $X$ and $f$ maps a source $w$
and a sink $v$ of $D$ to the same point $x_1\in X.$ Let
$D_{\sigma}$ be a graph obtained from $D$ by replacing
\begin{center}
\vspace*{.1in}
\parbox{.7in}{\psfig{figure=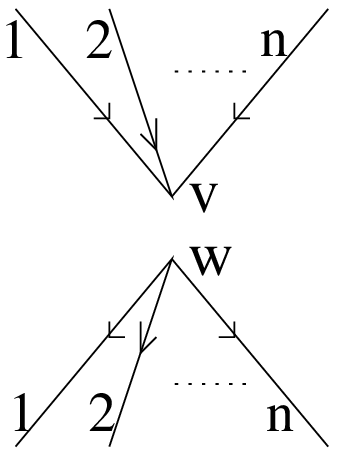,height=.7in}} by
\parbox{.7in}{\psfig{figure=diag1.eps,height=.7in}}.
\vspace*{.1in}
\end{center}
In the diagram we display $v$ and $w$ as separate points to
accentuate the fact that they are distinct in the domain of the
mapping $f.$
There is an obvious way to modify $f:D\to X$ to a map
$f_{\sigma}:D_{\sigma}\to X.$
We call $(f,\{f_{\sigma}\}_{\sigma\in S_n})$ a family of {\it skein related
graphs at $x_1$}.

\begin{definition}\label{2.2}\ \\
Let $X$ be a path connected topological space and let $I$ be the ideal
in $R{\cal G}_n(X)$ generated by two kinds of expressions:
\begin{enumerate}
\item $f-\sum_{\sigma\in S_n} \epsilon(\sigma) f_{\sigma},$ where
$\epsilon(\sigma)$ denotes the sign of $\sigma$ and
$(f,\{f_{\sigma}\}_{\sigma\in S_n})$ is a family of skein related
graphs at some point $x_1\in X.$
\item $L_e-n,$ where $e$ is the identity element in $\pi_1(X,x_0)$
(i.e. $L_e$ is homotopically trivial loop).
\end{enumerate}
Then the $R$-algebra ${\mathbb A}_n(X)=R{\cal G}_n(X)/I$ is called
the $n$-th skein algebra of $X.$
\end{definition}

Similarly we define ${\mathbb A}_n(X,x_0).$
Let $f_1:D_1\to X\times [0,1], f_2:D_2\to X\times [0,1]$ be elements of
${\cal G}_n(X,x_0).$ We define the product of them to be a map
$f_1\cdot f_2:D_1\cup D_2 \to X\times [0,1],$ such that
$$(f_1\cdot f_2)(d)=\begin{cases}
(x,{1\over 2}t) & \text{if $d\in D_1$ and $f_1(d)=
(x,t)$}\\ (x,{1\over 2}t+{1\over 2}) & \text{if $d\in D_2$ and
$f_2(d)=(x,t)$}\\
\end{cases}$$
This product extends to an associative (but generally
non-commutative) product in
$R{\cal G}_n(X,x_0).$ The identity in $R{\cal G}_n(X,x_0)$ is a
map $f:E\to X\times [0,1],$ where $E$ is a single edge and $f$ maps
$E$ onto $\{x_0\}\times [0,1].$

If $f:D\to X\times [0,1]$ is an element of ${\cal G}_n(X,x_0)$ such that
an $n$-valent source, $w,$ and an $n$-valent sink, $v,$ of $D$ are
mapped to a point $x_1\in X\times [0,1]$ then one can define
$D_{\sigma}$ and $f_{\sigma}: D_{\sigma}\to X\times [0,1],\ f_{\sigma}
\in {\cal G}_n(X,x_0),$ in exactly the same way as it was done for graphs in
${\cal G}_n(X)$ in the paragraph preceding Definition \ref{2.2}. We say,
as before, that $(f,\{f_{\sigma}\}_{\sigma\in S_n})$ are graphs skein related
at $x_1.$

\begin{definition}\label{2.3}\ \\
Let $X$ be a path connected topological space with a specified point
$x_0\in X,$ and let $I'$ be the ideal in $R{\cal G}_n(X,x_0)$ generated
by expressions
\begin{enumerate}
\item $f-\sum_{\sigma\in S_n} \epsilon(\sigma) f_{\sigma},$ where
$(f, \{f_{\sigma}\}_{\sigma\in S_n})$ is a family of skein related
graphs at some point $x_1\in X\times [0,1].$
\item $EL_e-n,$ where $e$ is the identity in $\pi_1(X,x_0).$
\end{enumerate}
Then the $R$-algebra ${\mathbb A}_n(X,x_0)=R{\cal G}_n(X,x_0)/I'$ is called
the $n$-th relative skein algebra of $(X,x_0).$
\end{definition}

Note that different choices of $x_0\in X$ give isomorphic
algebras ${\mathbb A}_n(X,x_0).$

Let $f\in {\cal G}_n(X),\ f:D\to X.$ Let $D'=D\cup E$ where $E$ is an
edge disjoint from $D.$
Then $D'$ has one $1$-valent sink $e_0$ and one
$1$-valent source $e_1,\ \{e_0,e_1\}=\partial E,$
and $D'\in {\cal G}_n'.$ We extend $f$ to $f':D'\to
X\times [0,1]$ in such a way that $f'(d)= (f(d),{1\over 2}),$ for
$d\in D,$ and
$f'(t)=(x_0,t),$ for $t\in [0,1]\simeq E,$ where $[0,1]\simeq E$ is a
preserving orientation parameterization of $E.$ This operation defines
an embedding $\imath:{\cal G}_n(X)\to {\cal G}_n(X,x_0),\ \imath(f)=f'$
Notice that $\imath$ induces a homomorphism
$\imath_*:{\mathbb A}_n(X)\to {\mathbb A}_n(X,x_0).$ Therefore we can
consider ${\mathbb A}_n(X,x_0)$ as an ${\mathbb A}_n(X)$-algebra.

Let $f:D\to X\times [0,1]$ be a map, $f\in {\cal G}_n(X,x_0).$
Let $\overline D\in {\cal G}_n$ be a graph obtained by identification of the
$1$-valent sink with the $1$-valent source in $D.$
Let us compose $f:D\to X\times [0,1]$ with a projection $X\times [0,1]\to X.$
This composition gives a map $\overline{f}:\overline{D}\to X,\ \overline{f}
\in {\cal G}_n(X).$
Therefore, we have a function $\overline{\cdot}: {\cal G}_n(X,x_0)\to
{\cal G}_n(X).$ This function can be extended to an $R$-linear homomorphism
${\mathbb T}: {\mathbb A}_n(X,x_0)\to {\mathbb A}_n(X).$ Notice that
for any graph $D\in {\cal G}_n(X),\ {\mathbb T}(\imath_*(D))$ is equal
to a union of $f:D\to X$ with a contractible loop in $X.$
Hence by Definition
\ref{2.2}(2) ${\mathbb T}(\imath_*(D))=n\cdot D.$ Since graphs in $X$ span
${\mathbb A}_n(X)$ the composition of $\imath_*:{\mathbb A}_n(X)\to
{\mathbb A}_n(X,x_0)$ with ${\mathbb T}: {\mathbb A}_n(X,x_0)\to {\mathbb
A}_n(X)$ is equal to $n$ times the identity on ${\mathbb A}_n(X).$ This
implies the following fact.

\begin{fact}\label{2.4}\ \\
If ${1\over n}\in R$ then $\imath_*:{\mathbb A}_n(X)\to {\mathbb
A}_n(X,x_0)$ is a monomorphism of rings.
\end{fact}

The next proposition summarizes basic properties of ${\mathbb
A}_n(X)$ and ${\mathbb A}_n(X,x_0).$

\begin{proposition}\label{2.5}\ \\
\vspace*{-0.2in}
\begin{enumerate}
\item The assignment $X\to {\mathbb A}_n(X)$ (respectively: $(X,x_0)\to
{\mathbb A}_n(X,x_0)$) defines a functor from the category of path
connected topological spaces (respectively: category of path connected
pointed spaces) to the category of commutative $R$-algebras
(respectively, the category of $R$-algebras).
\item If $f:X\to Y, f(x_0)=y_0,$ induces a surjection $f_*:
\pi_1(X,x_0)\to \pi_1(Y,y_0)$ then the corresponding
homomorphisms ${\mathbb A}_n(f):{\mathbb A}_n(X,x_0)\to {\mathbb
A}_n(Y,y_0),\\ {\mathbb A}_n(f):{\mathbb A}_n(X)\to {\mathbb A}_n(Y)$ are
epimorphisms of $R$-algebras.
\item The algebra ${\mathbb A}_n(X)$ is generated by loops in $X,$
i.e. by graphs $L_\gamma,$ for $\gamma\in \pi_1(X,x_0).$
\item The algebra ${\mathbb A}_n(X,x_0)$ is generated by graphs
$E_{g_i^{\pm 1}}$ and $EL_{\gamma},$ where $\{g_i\}_{i\in I}$
is a set of generators of $\pi_1(X,x_0)$ and $\gamma \in
\pi_1(X,x_0).$
\end{enumerate}
\end{proposition}

\begin{proof}
Since statements (1) and (2) of Proposition \ref{2.5} are obvious,
we give a proof of (3) and (4) only.

From the definition of a graph $D\in {\cal G}_n$ or $D\in{\cal G}_n'$ follows
that it has an equal number of $n$-valent sinks and sources.
Relation (1) of Definition \ref{2.2} and of Definition \ref{2.3} implies that
each pair of vertices of $f: D\to X,\ f\in {\cal G}_n(X)$
(respectively: of $f: D\to X\times [0,1],\ f\in {\cal G}_n(X,x_0)$) composed
of a sink and a source can be resolved and $f$ can be replaced by a linear
combination of graphs with a smaller number of sinks and sources. Therefore,
after a finite number of steps each graph in ${\cal G}_n(X)$ (respectively:
${\cal G}_n(X,x_0)$) can be expressed as a linear combination of graphs
without $n$-valent vertices.

(3) If $f:D\to X,\ f\in {\cal G}_n(X),$ and $D$ has no vertices then
$D$ is a union of loops, $D=S^1\cup S^1\cup ...\cup S^1,$ and therefore
$f=L_{\gamma_1}\cdot
L_{\gamma_2}\cdot...\cdot L_{\gamma_k}\in {\mathbb A}_n(X),$ for some
$\gamma_1,\gamma_2,...,\gamma_k\in \pi_1(X,x_0).$

(4) If $f:D\to X\times [0,1],\ f\in {\cal G}_n(X,x_0),$ and $D$ has no
$n$-valent vertices then $D=[0,1]\cup S^1\cup S^1\cup ...\cup S^1.$
Suppose that $[0,1]\stackrel{f}{\to} X\times [0,1]\to X$ represents
$\gamma_0\in \pi_1(X,x_0),$ and $f$ restricted to the $j$-th circle
represents the conjugacy
class of a $\gamma_j\in \pi_1(X,x_0),\ j=1,2,...,k.$ Then
$f=E_{\gamma_0}\cdot EL_{\gamma_1}\cdot EL_{\gamma_2}\cdot ...\cdot
EL_{\gamma_k}\in {\mathbb A}_n(X,x_0).$ Therefore ${\mathbb
A}_n(X,x_0)$ is generated by the elements $EL_{\gamma}$ and
$E_{\gamma'}, \ \gamma,\gamma'\in \pi_1(X,x_0).$ But each $E_{\gamma'}$ is
a product of elements $E_{g_i^{\pm 1}},$ where $\{g_i\}$ is a set of
generators of $\pi_1(X,x_0).$
\end{proof}

We will see later that ${\mathbb A}_n(X), {\mathbb A}_n(X,x_0)$ depend
only on $\pi_1(X,x_0).$
Moreover, if $\pi_1(X,x_0)$ is finitely generated group then the algebras
${\mathbb A}_n(X)$ and ${\mathbb A}_n(X,x_0)$ are also finitely generated.

Now we are ready to formulate the most important results of this paper.

\begin{theorem}\label{2.6}\ \\
Let $X$ be any (path connected) topological space and let $G=\pi_1(X,x_0),
x_0\in X.$ There are $R$-algebra homomorphisms
$$\Theta: {\mathbb A}_n(X,x_0)\to M_n(Rep_n^R(G))^{GL_n(R)},\
\theta: {\mathbb A}_n(X)\to Rep_n^R(G)^{GL_n(R)},$$ uniquely determined
by the following conditions:
\begin{enumerate}
\item $\Theta(E_{\gamma})=j_{G,n}(\gamma),\ \Theta(EL_{\gamma'})=
Tr(j_{G,n}(\gamma')),$ for any $\gamma,\gamma' \in \pi_1(X,x_0).$
\item $\theta(L_\gamma)=Tr(j_{G,n}(\gamma)),$
for any $\gamma\in \pi_1(X,x_0).$
\end{enumerate}

Moreover, the following diagram commutes:

\begin{equation}
\begin{tabular}{ccc}
${\mathbb A}_n(X,x_0)$ & $\stackrel{\Theta}{\longrightarrow}$ &
$M_n(Rep_n^R(G))^{GL_n(R)}$\\
$\Big \downarrow \scriptstyle{\mathbb T}$ & & $\Big \downarrow
\scriptstyle{Tr}$\\
${\mathbb A}_n(X)$ & $\stackrel{\theta}{\longrightarrow}$ &
$Rep_n^R(G)^{GL_n(R)}$\\
\label{e5}
\end{tabular}
\end{equation}
\end{theorem}

\begin{theorem}\label{2.7}\ \\
If $R$ is a field of characteristic $0$ then
$$\Theta:{\mathbb A}_n(X,x_0)\to M_n(Rep_n^R(G))^{GL_n(R)},\
\theta: {\mathbb A}_n(X)\to Rep_n^R(G)^{GL_n(R)},$$ are isomorphisms of
$R$-algebras.
\end{theorem}

Let $R$ be a field of characteristic $0.$ It can be shown that if
$X$ is a $3$-manifold then ${\mathbb A}_2(X)$ is isomorphic to the
Kauffman bracket skein module of $X,$ ${\cal S}_{2,\infty}(X,R,\pm
1).$ Moreover, if $X$ is a surface then ${\mathbb A}_2(X,x_0)$ is
isomorphic to the relative Kauffman bracket skein module of $X,$
${\cal S}^{rel}_{2,\infty}(X,R,\pm 1).$ See \cite{PS-2},
\cite{H-P}, for appropriate definitions and the notational
conventions. The main results of \cite{B-1}, \cite{B-2},
\cite{PS-1} and \cite{PS-2} relate the Kauffman bracket skein
modules of $3$-manifolds with the $SL_2$-representation theory of
their fundamental groups. Theorem \ref{2.7} generalizes these
results to groups $SL_n,$ for any $n.$

Moreover, it can be shown that in the case when $X$ is any path
connected topological space, ${\mathbb A}_2(X,x_0)$ and ${\mathbb A}_2(X)$
can be given the following simple algebraic description:
Let $G=\pi_1(X)$ and let $I$ be the ideal in the group ring $RG$
generated by elements
$h(g+g^{-1})-(g+g^{-1})h,$ where $g,h\in G.$ There is an involution
$\tau$ on $H(G)=RG/I$ sending $g$ to $g^{-1}.$ One can show
that ${\mathbb A}_2(X,x_0)$ is isomorphic
to $H(G)$ and ${\mathbb A}_2(X)$ is isomorphic to $H^+(G),$ where
$H^+(G)$ is the subring of $H(G)$
invariant under $\tau.$ The algebras $H(G), H^+(G)$ are introduced and
thoroughly investigated in \cite{B-H}. One of the main results of
\cite{B-H} is that
$H(G)=M_n(Rep_n^R(G))^{GL_n(R)}$ and $H^+(G)=Rep_n^R(G)^{GL_n(R)},$ for $n=2.$
(Compare also \cite{Sa-1}, \cite{Sa-2}). Theorem \ref{2.7} can be
considered as a generalization of this result to all values of $n.$

\section{Proof of Theorem \ref{2.6}}

Before we prove Theorem \ref{2.6} we give new definitions of
${\mathbb A}_n(X)$ and ${\mathbb A}_n(X,x_0)$ which only use
$G=\pi_1(X,x_0).$

Let $X$ be a path connected topological space and $x_0\in X.$
For any graph in ${\cal G}_n(X),$ i.e. a map $f:D\to X$ for some
$D\in {\cal G}_n,$ there is a map $f':D\to X$ homotopic to $f,$
which maps all vertices of $D$ to $x_0.$
Therefore the homotopy class of $f$ can be described by the graph
$D$ with each edge $E$ labeled by an element of $\pi_1(X,x_0)$ corresponding
to the map $f_{|E}': E\to X$ and each loop $L$ labeled by the
conjugacy class in $\pi_1(X,x_0)$ corresponding to the map
$f_{|L}':L\simeq S^1\stackrel{f}{\to} X.$ This description
does not need to be unique.

We denote the set of graphs in ${\cal G}_n$ all of whose edges are labeled by
elements of $G$ and all loops are labeled by conjugacy classes in $G$ by
${\cal G}_n(G).$ There is a natural multiplication operation on
${\cal G}_n(G).$ The product of $D_1,D_2\in {\cal G}_n(G)$ is the
disjoint union of $D_1$ and $D_2.$ Therefore $R{\cal G}_n(G)$ is a
commutative $R$-algebra with $\emptyset$ as the identity.

Let $D$ be a graph in ${\cal G}_n(G).$ We have noticed already that
$D$ corresponds to a map $f:D\to X$ which maps all vertices of
$D$ to $x_0$ and restricted to edges and loops of $D$ agrees with
their labeling. Such $f$ is unique up to a homotopy which fixes the
vertices of $D.$ Let $w$ be a source and $v$ be a sink in $D.$ Since
$f$ maps $v$ and $w$ to the same point in $X$ there exists a map
$f_{\sigma}:D_{\sigma}\to X$ defined for any $\sigma\in S_n$
as in the paragraph preceding Definition \ref{2.2}. Notice that
$f_{\sigma}$ maps all
vertices of $D_{\sigma}$ to $x_0\in X.$ Therefore, we can label all
edges of $D_{\sigma}$ by appropriate elements of $G$ and all loops of
$D_{\sigma}$ by appropriate conjugacy classes in $G,$ and hence consider
$D_{\sigma}$ as an element of ${\cal G}_n(G).$ Hence, we have showed
that one can replace any source $w$ and any sink $v$ in an arbitrary
graph $D\in {\cal G}_n(G)$
\begin{center}
\parbox{.7in}{\psfig{figure=diag3.eps,height=.7in}} by\quad
\parbox{.7in}{\psfig{figure=diag1.eps,height=.7in}}
\vspace*{.1in}
\end{center}
and obtain a well defined graph $D_{\sigma}\in {\cal G}_n(G).$

As an example consider the graph $D$ presented below.\\
\centerline{{\psfig{figure=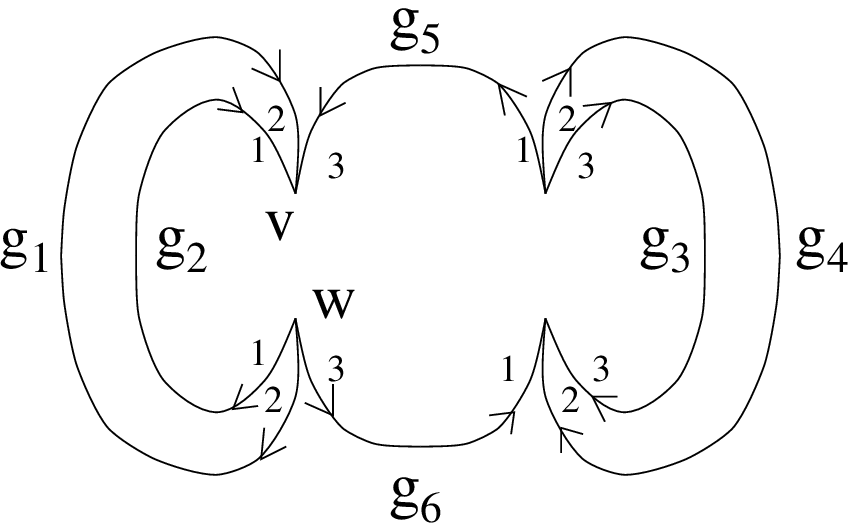,height=1.2in}}}\\
Replacing the vertices $v,w$ by a coupon decorated by $\sigma=(123)\in
S_3$ gives a diagram $D_{\sigma}:$\\
\centerline{\psfig{figure=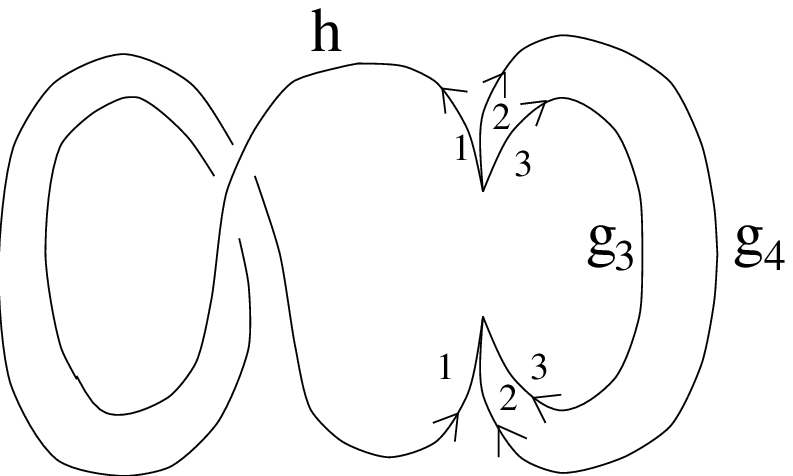,height=1in}}\\
where $h=g_6g_1g_2g_5.$

Now we are ready to define ${\mathbb A}_n(X)$ in terms of graphs in
${\cal G}_n(G).$ Namely, this algebra is isomorphic to
$R{\cal G}_n(\pi_1(X,x_0))/I,$ where $I\triangleleft R{\cal G}_n
(\pi_1(X,x_0))$ is an ideal generated by relations analogous to
relations (1) and (2) of Definition \ref{2.2} and by relations following from
the fact that the assignment ${\cal G}_n(G)\to {\cal G}_n(X)$ described above
is onto but not 1-1. The problem comes from the fact that one can take
a graph $D\in {\cal G}_n$ whose edges and loops are labeled in two
different ways such that the corresponding maps, $f, f':D\to X,$
sending the vertices of $D$ to $x_0$ are homotopic but
not by a homotopy relative to the vertices of $D.$ In
order to resolve this problem we need to allow an operation which
moves vertices of $D$ around paths in $X$ beginning and ending
at $x_0.$ Notice however that it suffices to move one vertex at the
time. The following fact summarizes our observations.

\begin{fact}\label{3.1}\ \\
Let $X$ be a path connected topological space with a
specified point $x_0\in X$ and let $G=\pi_1(X,x_0).$
Let $I$ be the ideal in $R{\cal G}_n(G)$ generated by
expressions of the following form:
\begin{equation} \parbox{.7in}{\psfig{figure=diag3.eps,height=.7in}}-
\sum_{\sigma\in S_n}\epsilon(\sigma)
\parbox{.7in}{\psfig{figure=diag1.eps,height=.7in}}. \label{e6}
\end{equation}
\begin{equation} \parbox{.4in}{\psfig{figure=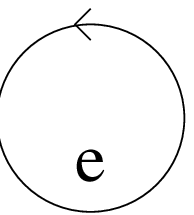,height=.4in}}-
n \label{e7}
\end{equation}
\begin{equation}
\vspace*{.2in} \psfig{figure=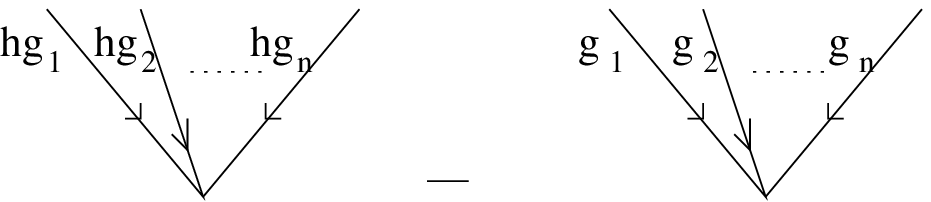,height=.5in} \label{e8}
\end{equation}
\begin{equation}
\vspace*{.2in} \psfig{figure=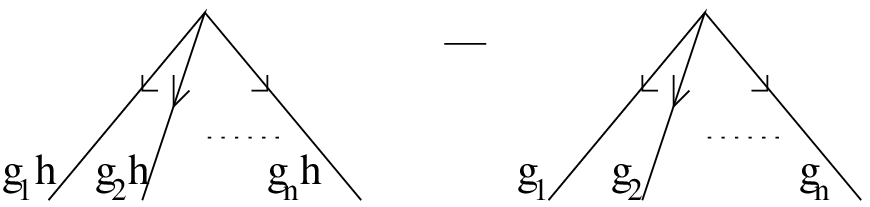,height=.5in} \label{e9}
\end{equation}
for any $g_1,g_2,...,g_n,h\in G.$

Then there is an isomorphism between $R$-algebras ${\mathbb A}_n(X)$ and
$R{\cal G}_n(G)/I$ assigning to each graph $f:D\to X, f\in {\cal
G}_n(X),$ with all vertices at $x_0$ the graph $D$ with every edge $E$
of $D$ decorated by the element of $\pi_1(X,x_0)$ corresponding to
$f_{|E}:E\to X,$ and every loop $L$ of $D$ decorated by the conjugacy class
in $\pi_1(X,x_0)$ represented by $f_{|L}:L\to X.$
\end{fact}

Now we will state a similar fact for ${\mathbb A}_n(X,x_0).$
Let ${\cal G}_n'(G)$ be a set of graphs in ${\cal G}_n'$ whose all edges are
labeled by elements of $G$ and whose all loops are labeled by conjugacy
classes in $G.$

There is a multiplication operation defined on ${\cal G}_n'(G)$ in the
following way. Let $D_1,D_2\in {\cal G}_n'(G),$ $v_i$ be the
$1$-valent source of $D_i,\ i\in\{1,2\},$ and let $w_i$ be the $1$-valent
sink of $D_i.$ Let $g_i$ be the label of the edge of $D_i$ joining
$v_i$ with $w_i.$
The graph $D_1\cdot D_2$ is obtained from the disjoint union of $D_1$ and
$D_2$ by identifying $v_1$ with $w_2.$ The edge of $D_1\cdot D_2$ joining
$v_2$ with $w_1$ is labeled by $g_1\cdot g_2.$ All other edges and loops of
$D_1\cdot D_2$ inherit labels form $D_1$ and $D_2.$ A single edge
labeled by $e\in G$ is the identity in ${\cal G}_n'(G).$

This multiplication extends to an associative (but not commutative)
multiplication in $R{\cal G}_n'(G).$

\begin{fact}\label{3.2}\ \\
Let $X$ be a path connected topological space with a specified point
$x_0\in X$ and let $G=\pi_1(X,x_0).$ Let $I'$ be the ideal in
$R{\cal G}_n'(G)$ generated by expressions (\ref{e6}),(\ref{e8}),
(\ref{e9}) and
$$\parbox{.4in}{\psfig{figure=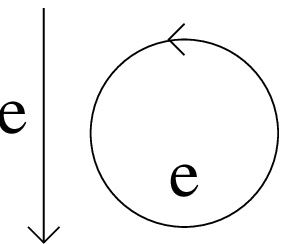,height=.4in}}-n\
\parbox{.4in}{\psfig{figure=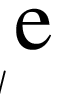,height=.4in}}.$$
Then ${\mathbb A}_n(X,x_0)\simeq R{\cal G}_n'(G)/I'.$
\end{fact}

Facts \ref{3.1} and \ref{3.2} show that the algebras ${\mathbb A}_n
(X,x_0),\ {\mathbb A}_n(X)$ depend only on $\pi_1(X,x_0).$ In fact
\ref{3.1} and \ref{3.2} give us models for ${\mathbb A}_n(X,x_0)$ and
${\mathbb A}_n(X)$ built from ${\cal G}_n(G)$ and ${\cal G}_n'(G).$ In
the rest of this section we will use these models.

Let us fix a commutative ring $R$ and a positive integer $n$ and a
topological space $X$ with $x_0\in X, \pi_1(X,x_0)=G.$
Let ${\cal R}=Rep_n^R(G)$ and let $V={\cal R}^n$ be a free $n$-dimensional
module over ${\cal R}$ with the standard basis, $\{e_1,e_2,...,e_n\},\
e_i=(0,0,...,1,...,0).$ The dual space $V^*$ has the dual basis
$e^1,e^2,...,e^n,\ e^i(e_j)=\delta_{i,j}.$ We will always use the
standard bases and therefore identify $V^*\otimes V\simeq End_{\cal
R}(V) \simeq M_n({\cal R}).$

Let $D$ be an element of ${\cal G}_n(G)$ or ${\cal G}_n'(G).$
We can decompose $D$ into arcs, sources and sinks:\\
\centerline{\psfig{figure=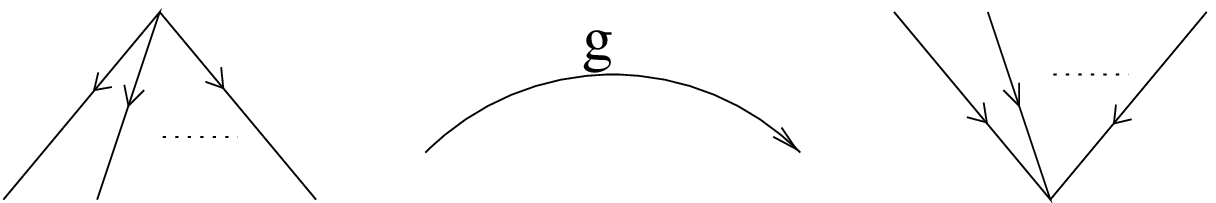,height=.5in}}

\begin{example}\label{3.3}\ \\
\parbox{0.8in}{\psfig{figure=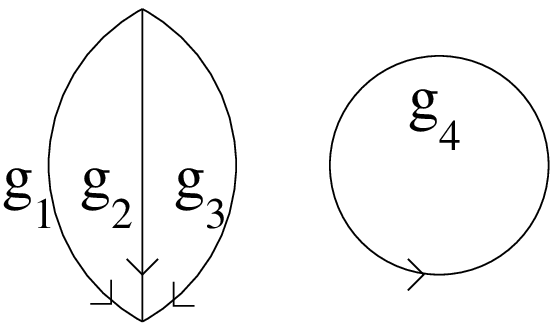,height=0.8in}} can be decomposed
to \ \parbox{0.9in}{\psfig{figure=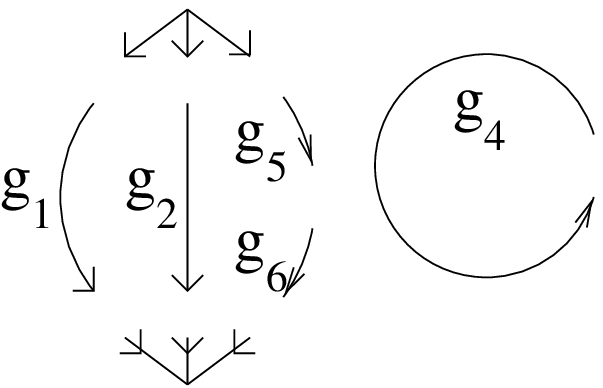,height=0.9in}}
where\\ $g_1,g_2,g_3,g_4,g_5,g_6\in G,\ g_3=g_6g_5.$
\end{example}

Notice that the decomposition of a graph is not unique since we allow
to cut each edge or loop into many pieces.

Let us assign to each $n$-valent source the tensor $\sum_{\sigma\in S_n}
\epsilon(\sigma) e_{\sigma(1)}\otimes e_{\sigma(2)}\otimes ...
\otimes e_{\sigma(n)}\in V^{\otimes n},$ to each $n$-valent sink the tensor
$\sum_{\sigma\in S_n} \epsilon(\sigma) e^{\sigma(1)}\otimes e^{\sigma(2)}
\otimes ... \otimes e^{\sigma(n)}\in (V^*)^{\otimes n}.$
To each edge labeled by $g$ we assign a tensor
in $V^*\otimes V\simeq End_{\cal R}(V)$ corresponding to
$j_{G,n}(g)\in SL_n({\cal R})\subset M_n({\cal R}).$

Let $D_0$ denote a graph $D$ decomposed into pieces.
We assign to $D_0$ the tensor product of
tensors corresponding to them. We denote this tensor by $T(D_0).$
Notice that $T(D_0)\in V^{\otimes N}\otimes (V^*)^{\otimes N},$ where
$N=$ the number of $1$-valent sources in $D_0\ =$ the number
of $1$-valent sinks in $D_0.$

Now we glue all components of $D_0$ together to get the graph $D$ back.
Whenever we glue an end of one piece with a beginning of another piece
in $D_0$ we make the corresponding contraction on $T(D_0).$
More specifically, suppose that the two free ends glued together
correspond to the two underlined components:
$$T(D_0)\in V\otimes ... \otimes \underline{V}\otimes ... \otimes V
\otimes V^*\otimes ... \otimes \underline{V^*}\otimes ... \otimes
V^*.$$ By applying to this tensor space the contraction map
$\underline{V}\otimes \underline{V^*}\to R$ (which is the
evaluation map $(v,f)\to f(v)$), we send $T(D_0)\in
V^{\otimes N}\otimes (V^*)^{\otimes N}$ to an element of
$V^{\otimes N-1}\otimes (V^*)^{\otimes N-1}.$
By repeating this process until we get the graph $D$ back,
we obtain an element of ${\cal R},$ if $D\in {\cal G}_n(G),$ or
an element of $M_n({\cal R}),$ if $D\in {\cal G}_n'(G).$
Notice that the above construction does not depend on the particular
decomposition of $D$ into pieces.
Therefore, we have defined functions:
\begin{equation}
\Theta: {\cal G}_n'(G)\to M_n({\cal R}),\quad  \theta: {\cal G}_n(G)\to
{\cal R}\label{e10}.
\end{equation}

\begin{lemma}\label{3.4}\ \\
Let $D$ be a graph in ${\cal G}_n(G)$ or in ${\cal G}_n'(G).$
Let $w$ be an $n$-valent source of $D$ and $v$ an $n$-valent sink of
$D.$ Let $D_{\sigma},$ for $\sigma\in S_n,$ be defined as at the beginning of
Section 4. If $D\in {\cal G}_n(G)$ then $\theta(D)=\sum_{\sigma\in S_n}
\epsilon(\sigma)\theta(D_{\sigma}).$ If $D\in {\cal G}_n'(G)$ then
$\Theta(D)=\sum_{\sigma\in S_n} \epsilon(\sigma)\Theta(D_{\sigma}).$
\end{lemma}

\begin{proof}
We will prove Lemma \ref{3.4} only for $D\in {\cal G}_n(G).$ For
$D\in {\cal G}_n'(G)$ the proof is identical.

Decompose $D$ and $D_{\sigma}$ into sources, sinks, and edges.
We denote the fragment of the decomposition of $D$, composed of the source $w$
and the sink $v$ by $D^0.$ We may assume that the decomposition of
$D_{\sigma}$ is identical to that of $D,$ except that it contains a
coupon $D^0_{\sigma}$ instead of $D^0.$
$$D^0=\parbox{.7in}{\psfig{figure=diag3.eps,height=.7in}}\hspace*{.7in}
D^0_{\sigma}= \parbox{.7in}{\psfig{figure=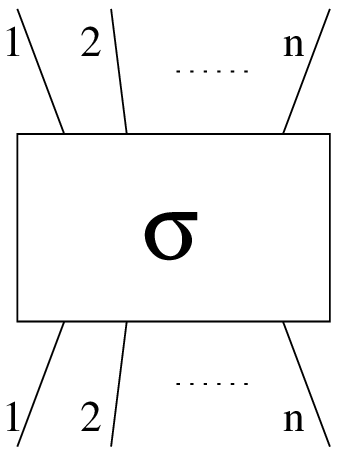,height=.7in}}$$
We order the $1$-valent sources and sinks of $D_{\sigma}^0$ consistently with
the ordering of the $1$-valent vertices of $D^0.$

Let $T(D^0)$ (respectively: $T(D_{\sigma}^0)$) be the tensor associated
to $D^0$ (respectively: $D^0_{\sigma}$). We assume that the $i$-th coordinate
of $T(D^0) \in V^{\otimes n}\otimes V^{*\otimes n}$ corresponds to the
$i$-th source of $D^0,$ if $1\leq i\leq n,$ or to $(i-n)$-th sink of
$D^0,$ if $n<i\leq 2n.$

Recall that $\theta(D), \theta(D_{\sigma})\in {\cal R}$ are results of
contractions of tensors associated with elements of decompositions of
$D$ and $D_{\sigma}.$ Since the decompositions of $D$ and $D_{\sigma}$
chosen by us differ only by elements $D^0, D^0_{\sigma},$ the proof
of Lemma \ref{3.4}
can be reduced to a local computation on tensors. Namely, it is enough
to prove that
\begin{equation}
T(D^0)=\sum_{\sigma\in S_n}\epsilon(\sigma)T(D^0_{\sigma}) \label{e11}
\end{equation}
Notice that each edge of $D^0_{\sigma}$ is labeled by the identity map in
$End_R(V).$ This map is represented by $\sum_{i=1}^n e_i\otimes e^i$ in
$V\otimes V^*\simeq End_R(V).$ Therefore, if $\sigma=id\in S_n$ then
$$T(D^0_{\sigma})=\sum_{i_1,i_2,...,i_n\in\{1,2,...,n\}}
e_{i_1}\otimes e_{i_2}\otimes ... \otimes e_{i_n}\otimes
e^{i_1}\otimes e^{i_2}\otimes ... \otimes e^{i_n}.$$
Similarly, for any $\sigma\in S_n,$ we have
$$T(D^0_{\sigma})=\sum_{i_1,i_2,...,i_n\in\{1,2,...,n\}}
e_{i_1}\otimes e_{i_2}\otimes ... \otimes e_{i_n}\otimes
e^{i_{\sigma(1)}}\otimes e^{i_{\sigma(2)}}\otimes ... \otimes
e^{i_{\sigma(n)}}.$$
Therefore,
\begin{equation}
\sum_{\sigma\in S_n}\epsilon(\sigma)T(D^0_{\sigma}) =
\hspace*{-.2in}\sum_{\sigma\in S_n\atop i_1,i_2,...,i_n\in\{1,2,...,n\}}
\hspace*{-.2in}\epsilon(\sigma) e_{i_1}\otimes e_{i_2}\otimes ...
\otimes e_{i_n}\otimes
e^{i_{\sigma(1)}}\otimes e^{i_{\sigma(2)}}\otimes ... \otimes
e^{i_{\sigma(n)}}\label{e12}
\end{equation}
Note that we can assume that the numbers $i_1,i_2,...,i_n$ appearing
on the right side
of (\ref{e12}) are all different. Indeed, if $i_j=i_k, j\ne k,$ then
there is an equal number of even and odd permutations contributing the
term $$e_{i_1}\otimes e_{i_2}\otimes ...\otimes e_{i_n}\otimes
e^{j_1}\otimes e^{j_2}\otimes ... \otimes e^{j_n}$$ to the sum on the
right side of (\ref{e12}), for any $j_1, j_2,..., j_n.$

%

Therefore, we can assume that the numbers $(i_1,i_2,...,i_n)$ appearing in
each term of the sum on the right side of (\ref{e12}) form a
permutation $\tau$ of $(1,2,...,n).$ Hence we have,
$$\sum_{\sigma\in S_n}\epsilon(\sigma)T(D^0_{\sigma})=
\sum_{\sigma,\tau \in S_n}
\epsilon(\sigma) e_{\tau(1)}\otimes e_{\tau(2)}\otimes ... \otimes e_{\tau(n)}
\otimes e^{\tau(\sigma(1))}\otimes e^{\tau(\sigma(2))}\otimes ... \otimes
e^{\tau(\sigma(n))}.$$

Substitute $\tau'$ for $\tau\circ \sigma.$ Then $\epsilon(\sigma)=
\epsilon(\tau)\epsilon(\tau')$ and we get
$$\sum_{\sigma\in S_n}\epsilon(\sigma)T(D^0_{\sigma})=
\sum_{\tau,\tau' \in S_n}
\epsilon(\tau)\epsilon(\tau') e_{\tau(1)}\otimes e_{\tau(2)}\otimes ...
\otimes e_{\tau(n)} \otimes e^{\tau'(1)}\otimes e^{\tau'(2)}\otimes ...
\otimes e^{\tau'(n)}.$$
Notice that the right side of the above equation is equal to
$$\left( \sum_{\tau\in S_n} \epsilon(\tau) e_{\tau(1)}\otimes e_{\tau(2)}
\otimes ... \otimes e_{\tau(n)}\right) \otimes
\left( \sum_{\tau'\in S_n} \epsilon(\tau') e^{\tau'(1)}\otimes e^{\tau'(2)}
\otimes ... \otimes e^{\tau'(n)}\right).$$
But the expression above is exactly the tensor assigned to\\
\centerline{\psfig{figure=diag3.eps,height=.8in}}

Therefore we have proved (\ref{e11}) and completed the proof of Lemma
\ref{3.4} \end{proof}

The study of $SL_n$-actions on linear spaces was one of the main
objectives of classical invariant theory. In particular, Weyl
(\cite{We}) determined all invariants of the action of $SL(V)$ on
$V\otimes ...\otimes V\otimes V^*\otimes ...\otimes V^*$ and gave a
full description of relations between them.
The set of ``typical'' invariants consists of
brackets $[v_1,...,v_n]=Det(v_1,...,v_n), [\phi_1,...,\phi_n]^*=
Det(\phi_1,...,\phi_n),$ where $v_1,...,v_n\in V,
\phi_1,...,\phi_n\in V^*,$ and contractions $\phi_j(v_i).$
The identity
\begin{equation}
[v_1,...,v_n][\phi_1,...,\phi_n]^*=Det(\phi_j(v_i))_{i,j=1}^n
\label{e13}
\end{equation}
is one of the fundamental identities relating the typical invariants.
The bracket $[\cdot,\cdot,...,\cdot]$ is a skew symmetric linear functional on
$V\otimes V\otimes ...\otimes V$ and hence an element of $\bigwedge^n
V^*.$ Similarly, $[\cdot,\cdot,...,\cdot]^*\in \bigwedge^n V.$
Note that the sources and sinks of graphs considered by us are labeled
exactly by the tensors $[\cdot,\cdot,...,\cdot]^*$ and
$[\cdot,\cdot,...,\cdot].$ (However, $V$ is in our case a free module
over ${\cal R}=Rep_n^R(G)$.)

It follows from the proof of Lemma \ref{3.4} that
$$\sum_{\sigma\in S_n} \epsilon(\sigma)\ \parbox{.7in}
{\psfig{figure=diag1.eps,height=.7in}}$$
represents the tensor in $Hom(V\otimes...\otimes V\otimes
V^*\otimes...\otimes V^*, {\cal R})=V^*\otimes...\otimes V^*\otimes
V\otimes...\otimes V$ assigning to $(v_1,v_2,...,v_n,\phi_1,\phi_2,
...\phi_n)$ the value $Det(\phi_j(v_i))_{i,j=1}^n.$ Therefore, the
identity $$\theta(D)=\sum_{\sigma\in S_n}
\epsilon(\sigma)\theta(D_{\sigma})$$
is essentially equivalent to (\ref{e13}).

\begin{lemma}\label{3.5}\ \\
Let $L_g,E_g,EL_g$ be graphs defined as in Section 3 but considered as
elements of ${\cal G}_n(G)$ and ${\cal G}_n'(G)$ i.e.
\begin{enumerate}
\item $L_g\in {\cal G}_n(G)$ is a single loop labeled by the conjugacy class
of $g\in G,$
\item $E_g\in {\cal G}_n'(G)$ is a single edge labeled by $g\in G,$
\item $EL_g\in {\cal G}_n'(G)$ is a graph composed of an
edge labeled by the identity in $G$ and of a loop labeled by the
conjugacy class of $g\in G.$
\end{enumerate}
Under the above assumptions the functions $\Theta, \theta$ satisfy
conditions (1) and (2) of Theorem \ref{2.6}.
\end{lemma}

\begin{proof}
\begin{enumerate}
\item $L_g$ can be decomposed into a single arc\ \parbox{.4in}
{\psfig{figure=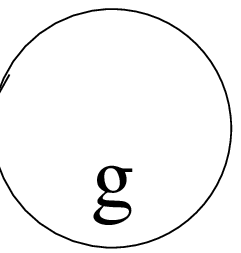,height=.4in}}
which has associated the tensor $j_{G,n}(g)\in SL_n({\cal R})\subset
V^*\otimes V.$ The contraction of this tensor gives
$\theta(L_g)=Tr(j_{G,n}(g)).$
\item $E_g$ is a single arc. Therefore
$\Theta(E_g)=T(E_g)=j_{G,n}(g).$
\item $EL_g$ can be decomposed into
\parbox{.4in}{\psfig{figure=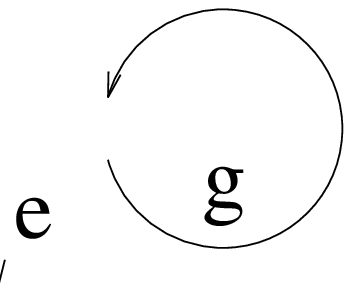,height=.4in}}.
The tensor associated with this decomposition is
$id\otimes j_{G,n}(g)\in End(V)\otimes End(V).$ After making a
contraction corresponding to the identification of the ends of the arc
we get
$$\Theta(EL_g)=id\cdot Tr(j_{G,n}(g))\in End(V).$$
\end{enumerate}
\end{proof}

\begin{lemma}\label{3.6}\ \\
Let $D,D'\in {\cal G}_n(G)$ be two graphs which are identical as
unlabeled graphs and which have the same labeling of edges and loops
except the labeling of edges incident to a vertex $v.$ Moreover,
suppose that
\begin{enumerate}
\item if $v$ is a source then the edges in $D$ incident to $v$ are labeled by
$g_1,g_2,...,g_n$ and the edges in $D'$ incident to $v$ are labeled
by $g_1h,g_2h,...,g_nh$ for some $g_1,g_2,...,g_n,h\in G.$
\item if $v$ is a sink then the edges in $D$ incident to $v$ are labeled by
$g_1,g_2,...,g_n$ and the edges in $D'$ incident to $v$ are labeled by
$hg_1,hg_2,...,hg_n$ for some $g_1,g_2,...,g_n,h\in G.$
\end{enumerate}
Under the above assumptions $\theta(D)=\theta(D').$ An analogous fact
is true for graphs in ${\cal G}_n'(G).$
\end{lemma}

\begin{proof}
We prove part (1) only. The proof of part (2) is analogous.

Let $v$ be a source.
Notice that $D$ and $D'$ have identical decompositions into sinks,
sources, and arcs except that $D_0=$\parbox{.3in}
{\psfig{figure=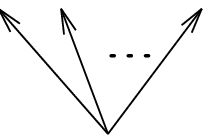,height=.3in}} is an element of the
decomposition of $D$ and the diagram $D_0'=$\parbox{.6in}
{\psfig{figure=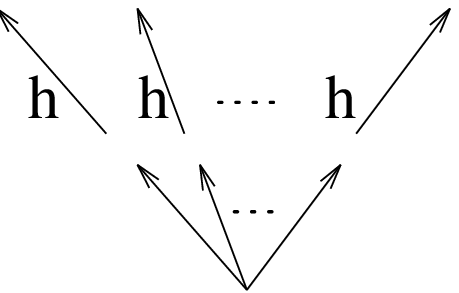,height=.6in}} is a fragment of a
decomposition of $D'.$ Therefore we need to show that the tensors
assigned to the above diagrams are identical. Notice that the tensor
associated to $D_0,\ T(D_0),$ is an element of the one dimensional
${\cal R}$-linear space of skew symmetric tensors $\bigwedge^n
V\subset V^n.$ Let $A:V\to V$ be an endomorphism given in the standard
coordinates of $V$ by $j_{G,n}(h)\in SL_n({\cal R}).$ $A$ induces an
endomorphism $\wedge A:\bigwedge^n V\to \bigwedge^n V$ with the
property that $\wedge A(T(D_0))=Det(A)T(D_0)\in \bigwedge^n V.$
Notice that $\wedge A(T(D_0))$ is exactly the tensor associated to
$D_0'.$ Since $Det(A)=1,$ the tensors associated to $D_0$ and
$D_0'$ are equal.
\end{proof}

{\bf Proof of Theorem \ref{2.6}}\\
Let us extend the functions $\theta$ and $\Theta$ to all $R$-linear
combinations of graphs in ${\cal G}_n(G)$ and ${\cal G}_n'(G)$
respectively. Facts \ref{3.1} and \ref{3.2} and Lemmas \ref{3.4},
\ref{3.5}, and \ref{3.6} imply that these functions descend to
$R$-linear homomorphisms
$$\Theta: {\mathbb A}_n(X,x_0)\to M_n({\cal R}),\qquad
\theta: {\mathbb A}_n(X)\to {\cal R}.$$
By Lemma \ref{3.5}, $\Theta$ and $\theta$ satisfy the conditions (1)
and (2) of Theorem \ref{2.6}.

We have showed in proposition \ref{2.5}(4) that ${\mathbb A}_n(X,x_0)$
is generated by elements $E_{\gamma}$ and $EL_{\gamma'},$ for
$\gamma, \gamma'\in G=\pi_1(X,x_0).$ By Proposition \ref{1.3} and the
paragraph preceding it, $\Theta(E_\gamma)=j_{G,n}(\gamma)$ and
$\Theta(EL_{\gamma'})=Tr(j_{G,n}(\gamma'))$
belong to $M_n(Rep_n^R(G))^{GL_n(R)}.$ Therefore the image of $\Theta$
lies in $M_n(Rep_n^R(G))^{GL_n(R)}.$ We show analogously that the image
of $\theta$ lies in $Rep_n^R(G)^{GL_n(R)}.$ Therefore the proof will be
completed if we show that the diagram of Theorem \ref{2.6} commutes.

Let $D\in {\cal G}_n'(G),$ $G=\pi_1(X,x_0),$ represent an
element of ${\mathbb A}_n(X,x_0).$ Then $\Theta(D)\in
M_n(Rep_n^R(G))^{GL_n(R)}$ is the result of a contraction of tensors
associated with elements of some decomposition of $D.$ Notice that
${\mathbb T}(D)$ is an element of ${\mathbb A}_n(X)$ represented by
the diagram $D$ with its $1$-valent vertices identified
\footnote{Recall that ${\mathbb T}$ was defined in the paragraph
preceding Fact \ref{2.4}.}. Hence $\theta({\mathbb T}(D))$ is a
contraction of $\Theta(D),$ i.e. $\theta({\mathbb
T}(D))=Tr(\Theta(D)).$ Since the elements $D\in {\cal G}_n'(G)$
span ${\mathbb A}_n(X,x_0)$ the diagram of Theorem \ref{2.6}
commutes.
\boxx

\section{Proof of Theorem \ref{2.7}}

Since now we assume that $R$ is a field of characteristic $0.$

At the beginning we state the first and the second fundamental theorems
of invariant theory following the approach of Procesi, \cite{Pro-1}
(Compare also \cite{Ra}).

Let $I$ be an infinite set.
Let $P_n(I)$ and $A_i$ be as before,
$$P_n(I)=R[x^i_{jk},\ j,k\in\{1,2,...,n\}, i\in I],\ A_i=(x^i_{jk})\in
M_n(P_n(I)).$$
We are going to present Procesi's description of the ring
$M_n(P_n(I))^{GL_n(R)}.$

Let $T$ be a commutative $R$-algebra freely generated by the
symbols $Tr(X_{i_1}X_{i_2}...X_{i_k}),$ where $i_1,i_2,...,i_k\in I.$
We adopt the convention that $Tr(M)=Tr(N)$ if and only if the monomial $N$ is
obtained from $M$ by a cyclic permutation.
Let $T\{X_i\}_{i\in I}$ be a non-commutative $T$-algebra freely generated by
variables $X_i,\ i\in I.$ We have a natural $T$-linear homomorphism
$Tr:T\{X_i\}_{i\in I}\to T$ which assigns to $X_{i_1}X_{i_2}...X_{i_k}\in
T\{X_i\}_{i\in I}$ an element $Tr(X_{i_1}X_{i_2}...X_{i_k})\in T.$

There is a homomorphism $\pi:T\{X_i\}_{i\in I}\to M_n(P_n(I))$ uniquely
determined by the conditions:
\begin{itemize}
\item $\pi(X_i)=A_i$
\item $\pi(Tr(X_{i_1}X_{i_2}...X_{i_k}))=Tr(A_{i_1}A_{i_2}...A_{i_k})\in
P_n(I)\subset M_n(P_n(I)).$\footnote{We identify $P_n(I)$ with the
scalar matrices in $M_n(P_n(I)).$}
\end{itemize}
Proposition \ref{1.2} and Lemma \ref{1.1} imply that the image of
$\pi$ is fixed by the $GL_n(R)$-action on $M_n(P_n(I)),$ i.e.
$\pi:T\{X_i\}_{i\in I}\to M_n(P_n(I))^{GL_n(R)}.$

Notice that the following digram commutes.

\begin{center}
\begin{tabular}{ccc}
$T\{X_i\}_{i\in I}$ & $\stackrel{\pi}{\longrightarrow}$ &
$M_n(P_n(I))^{GL_n(R)}
$\\ $\Big\downarrow \scriptstyle{Tr}$ & & $\Big\downarrow \scriptstyle{Tr}$\\
$\hspace*{-.1in} T$ & $\stackrel{\pi_{|T}}{\longrightarrow}$ &
$P_n(I)^{GL_n(R)}$\\
\end{tabular}
\end{center}

The following version of {\it The First Fundamental Theorem} of
invariant theory of $n\times n$ matrices is due to Procesi, \cite{Pro-1}.

\begin{theorem}\label{4.1}\ \\
$\pi:T\{X_i\}_{i\in I}\to M_n(P_n(I))^{GL_n(R)}$ is an epimorphism.
\end{theorem}

Before we state the second fundamental theorem of invariant theory of
$n\times n$ matrices we need some preparations.

Suppose that $\{1,2,...,m\}\subset I$ and specify $i_0\in \{1,2,...,m\}.$
We can present any $\sigma\in S_m$ as a product of cycles in such a way that
$i_0$ is the first element of the first cycle,
$\sigma=(i_0,i_1,...,i_s)(j_0,j_1,...,j_t)...(k_0,k_1,...,k_v).$
We define $\Phi_{\sigma,i_0}(X_1,X_2,...,X_m)$ to be equal to
$$X_{i_0}X_{i_1}...X_{i_s}Tr(X_{j_0}X_{j_1}...
X_{j_t})...Tr(X_{k_0}X_{k_1}...X_{k_v})\in T\{X_i\}_{i\in I}.$$
We also define another expression which does not depend on $i_0,$
$$\Phi_{\sigma}(X_1,X_2,...,X_m)=Tr(X_{i_0}X_{i_1}...X_{i_s})Tr(X_{j_0}X_{j_1}...X_{j_t})...Tr(X_{k_0}X_{k_1}...X_{k_v})\in T.$$

Let $F(X_1,X_2,...,X_m)=\sum_{\sigma\in S_m} \epsilon(\sigma)
\Phi_{\sigma}(X_1,X_2,...,X_m)\in T.$\\
$F(X_1,X_2,...,X_{n+1})$ is called the fundamental trace identity of
$n\times n$ matrices.

Procesi argues that there exists a unique element
$G(X_1,X_2,...,X_n)\in T\{X_i\}_{i\in I}$ involving only the variables
$X_1,...,X_n$ and the traces of monomials in these variables such that
$$F(X_1,X_2,...,X_{n+1})=Tr(G(X_1,X_2,...,X_n)X_{n+1})\in T\{X_i\}_{i\in I}.$$
Procesi gives an explicit formula for $G(X_1,X_2,...,X_n),$ but we
want to give a different formula, which will be more suitable for our
purposes.

\begin{lemma}\label{4.2}\ \\
$$G(X_1,X_2,...,X_n)=
\sum_{\sigma\in S_n}\epsilon(\sigma)\Phi_{\sigma}(X_1,X_2,...,X_n)-
\hspace*{-.1in}\sum_{i\in \{1,2,...,n\}\atop \sigma\in S_n}
\hspace*{-.2in}\epsilon(\sigma)
\Phi_{\sigma,i}(X_1,X_2,...,X_n).$$
\end{lemma}

\begin{proof}
It follows from the remarks preceding Lemma \ref{4.2} that it is enough to
show that if we multiply the right side of the equation of
Lemma \ref{4.2} by $X_{n+1}$ then the trace of it will be equal to
$F(X_1,...,X_{n+1}),$ i.e. we have to prove that
$$\sum_{\sigma\in S_n}\epsilon(\sigma)\Phi_{\sigma}(X_1,X_2,...,X_n)
Tr(X_{n+1})-\sum_{i\in \{1,2,...,n\}\atop \sigma\in S_n}
\hspace*{-.2in} \epsilon(\sigma)
Tr(\Phi_{\sigma,i}(X_1,X_2,...,X_n)X_{n+1})$$
\begin{equation}
 =F(X_1,X_2,...,X_{n+1}).\label{e14}
\end{equation}

Notice that $\epsilon(\sigma)\Phi_{\sigma}(X_1,X_2,...,X_n)Tr(X_{n+1})=
\epsilon(\sigma')\Phi_{\sigma'}(X_1,X_2,...,X_n,X_{n+1}),$ where
$\sigma'\in S_{n+1},\ \sigma'(i)=\sigma(i),$ for $i\in
\{1,2,...,n\},$ and $\sigma'(n+1)=n+1.$

Similarly we can simplify $Tr(\Phi_{\sigma,i}(X_1,X_2,...,X_n)X_{n+1}).$
Suppose that $$\sigma=(i_0,i_1,...,i_s)(j_0,j_1,...,j_t)...(k_0,k_1,...,k_v)
\in S_n,$$ where $i_0=i.$ Then
$$Tr(\Phi_{\sigma,i}(X_1,X_2,...,X_n)X_{n+1})= \Phi_{\sigma'}
(X_1,X_2,...,X_n,X_{n+1}),$$ for
$\sigma'=(i_0,i_1,...,i_s,n+1)(j_0,j_1,...,j_t)...(k_0,k_1,...,k_v)
\in S_{n+1}.$ Notice that $\epsilon(\sigma')=-\epsilon(\sigma).$ Therefore
the left side of the equation (\ref{e14}) is equal to
$$\hspace*{-.2in} \sum_{\sigma'\in S_{n+1},\atop \sigma'(n+1)=n+1}
\hspace*{-.2in} \epsilon(\sigma')\Phi_{\sigma}
(X_1,X_2,...,X_n,X_{n+1})+ \hspace*{-.2in} \sum_{i\in\{1,2,...,n\},
\ \sigma'\in S_{n+1}\atop {\rm such\ that\ } \sigma'(n+1)=i}
\hspace*{-.3in} \epsilon(\sigma')
\Phi_{\sigma}(X_1,X_2,...,X_n,X_{n+1}).$$
The above expression is obviously equal to $F(X_1,X_2,...,X_{n+1}).$
\end{proof}

Now we are ready to state {\it The Second Fundamental Theorem} of
invariant theory of $n\times n$ matrices, \cite{Pro-1}.

\begin{theorem}\label{4.3}\ \\
The kernel of $\pi$ is generated by elements
$G(M_1,M_2,...,M_n),F(N_1,N_2,...,N_{n+1}),$ where $M_1,M_2,...,M_n,
N_1,N_2,...,N_{n+1}$ are all possible monomials in the variables
$X_i,\ i\in I.$
\end{theorem}

Let $X$ be a path connected topological space.
We choose a presentation $<g_i\ i\in I|r_j\ j\in J>$ of $G=\pi_1(X,x_0)$
such that
\begin{itemize}
\item $I$ is an infinite set
\item the inverse of every generator is also a generator.
\item the defining relations $r_j$ are products of
non-negative powers of generators.
\end{itemize}
Note that such presentation always exists (even if $G$ is finitely
generated).

Let $\psi: T\{X_i\}_{i\in I}\to {\mathbb A}_n(X,x_0)$ be an $R$-homomorphism
such that $\psi(X_i)=E_{g_i}$ and $\psi(Tr(X_{i_1}X_{i_2}...X_{i_k}))=
EL_{g_{i_1}g_{i_2}...g_{i_k}}.$ Recall that by Proposition \ref{2.4},
${\mathbb A}_n(X)$ can be considered as a subalgebra of ${\mathbb
A}_n(X,x_0),$ in such a way that $L_{\gamma}\in {\mathbb A}_n(X)$
is identified with $EL_{\gamma}\in {\mathbb A}_n(X,x_0).$
Hence
$\psi(Tr(X_{i_1}X_{i_2}...X_{i_k}))\in {\mathbb A}_n(X)$ and
$\psi$ restricts to $\psi:T\to {\mathbb A}_n(X).$ Moreover the
following diagram commutes.
\begin{equation}
\begin{tabular}{ccc}
$T\{X_i\}_{i\in I}$ & $\stackrel{\psi}{\longrightarrow}$ &
${\mathbb A}_n(X,x_0)$\\
$\Big\downarrow \scriptstyle{Tr}$ & & $\Big\downarrow
\scriptstyle{Tr}$\\
\hspace*{-.2in}$T$ & $\stackrel{\psi}{\longrightarrow}$ & ${\mathbb A}_n(X)$\\
\end{tabular} \label{e15}
\end{equation}
We are going to show that the kernel of $\psi:T\{X_i\}_{i\in I} \to
{\mathbb A}_n(X,x_0)$ contains the kernel of $\pi: T\{X_i\}_{i\in I}\to
M_n(P_n(I))^{GL_n(R)}$ and therefore $\psi$ descends to
a homomorphism $M_n(P_n(I))^{GL_n(R)}\to {\mathbb A}_n(X,x_0).$

We will need the following fact due to Formanek (Proposition 45
\cite{For}).

\begin{proposition}\label{4.4}\ \\
For any matrix $A\in M_n(R)$
$$Det(A)={1\over n!}\sum_{\sigma\in S_n} \epsilon(\sigma)Tr(A^{c_1})
Tr(A^{c_2})...Tr(A^{c_k}),$$ where $c_1,c_2,...,c_k$ denote the lengths of
all cycles in $\sigma.$
\end{proposition}

For completeness we sketch a proof of Proposition \ref{4.4}.
A multilinearization of
the determinant, $Det: M_n(R)\to R,$ gives a function on $n$-tuples
of $n\times n$ matrices
$${\cal M}(X_1,...,X_n)=\sum_{\sigma\in S_n} Det(X_{\sigma}),$$ where
$X_{\sigma}$ is a matrix whose $i$-th row is the $i$-th row of
$X_{\sigma(i)}.$ Note that ${\cal M}(A,...,A)=n! Det(A)$ and
therefore the identity of Proposition \ref{4.4} is a special case of
the following identity:
$$ {\cal M}(X_1,...,X_n)=\sum_{\sigma\in S_n} \epsilon(\sigma)
\Phi_{\sigma}(X_1,...,X_n),$$ where $\Phi_{\sigma}$ was
defined in the second paragraph after Theorem \ref{4.1}. Formanek
gives the following proof of the above identity.
Assume that $1,2,...,n\in I.$ Since $A_1,...,A_n\in M_n(P_n(I))$
represent generic matrices, in order to prove the above identity it is
enough to show it for $X_1=A_1, ..., X_n=A_n.$
Since ${\cal M}(A_1,...,A_n)$ is an invariant polynomial function on
$n$-tuples of matrices, the First Fundamental Theorem implies that
${\cal M}(A_1,...,A_n)$ can be expressed in terms of
traces of monomials in $A_1,...,A_n.$ Since ${\cal
M}(A_1,...,A_n)$ is linear with respect to $A_1,...,A_n,$ it
is a linear combination of terms $Tr(A_{i_1}...A_{i_s})...Tr(A_{j_1}...
A_{j_t}),$ where $i_1,...i_s,...,j_1,...,j_t$ form a permutation of
$1,2,...,n.$ Therefore \\
${\cal M}(A_1,...,A_n)=\sum_{\sigma\in S_n} \alpha_{\sigma}
\Phi_{\sigma}(A_1,...,A_n)$ and hence
\begin{equation}
{\cal M}(X_1,...,X_n)=\sum_{\sigma\in S_n} \alpha_{\sigma}
\Phi_{\sigma}(X_1,...,X_n), \label{e15a}
\end{equation}
for any $n\times n$ matrices $X_1,...,X_n.$
We need to prove that $\alpha_{\sigma}=\epsilon(\sigma).$
If we restrict the above equation to matrices $A_1,..,A_n\in
M_{n-1}(P_{n-1}(I))$
embedded into $M_n(P_{n-1}(I))$ in the standard, non unit preserving way, we
will get the following polynomial identity on $(n-1)\times (n-1)$ matrices:
$$\sum_{\sigma\in S_n} \alpha_{\sigma} \Phi(A_1,...,A_n)=0.$$
It is not difficult to see that the Second Fundamental Theorem implies
that $F(A_1,...,A_n)$ is the only (up to scalar) $n$-linear trace
identity of degree $n$ on matrices $A_1,...A_n\in M_{n-1}(P_{n-1}(I)).$
Therefore $\alpha_{\sigma}=\epsilon(\sigma)c,$ for some fixed $c.$
Substituting the matrix $(x_{ij})$ with a single nonzero entry
$x_{ii}=1$ for $X_i$ in (\ref{e15a}) we get $c=1.$ Thus the proof of
Proposition \ref{4.4} is finished.

The specialization $A=Id\in M_n(R)$ in Proposition \ref{4.4} yields
the following corollary.

\begin{corollary}\label{4.5}\ \\
For any positive integer $n$
$$\sum_{\sigma\in S_n} \epsilon(\sigma) n^{c(\sigma)}=n!,$$
where $c(\sigma)$ is the number of cycles in the cycle
decomposition of $\sigma.$
\end{corollary}

From the definition of $T\{X_i\}_{i\in I}$ it immediately follows that for any
family of matrices $\{M_i\}_{i\in I}\subset M_n(R),$ there is well defined
substitution
$$X_i\to M_i,\quad \quad Tr(X_{i_1}X_{i_2}...X_{i_k})\to
Tr(M_{i_1}M_{i_2}...M_{i_k})\in R\subset M_n(R),$$
which can be extended to the whole ring $T\{X_i\}_{i\in I}.$
Therefore, if $H(X_{i_1},X_{i_2}, ..., X_{i_k})$ is an element of
$T\{X_i\}_{i\in I}$ involving variables $X_{i_1},X_{i_2}, ..., X_{i_k}$ then
$H(M_{i_1},M_{i_2}, ..., M_{i_k})$ is a well defined  matrix in $M_n(R).$

\begin{lemma}\label{4.6}\ \\
If $N_1,N_2,...,N_n$ are any monomials in the variables $X_i,\ i\in I,$
then\\ $\psi(G(N_1,N_2,...,N_n))=0.$
\end{lemma}

\begin{proof}
By the definition of $\psi$ (given in the second paragraph after
Theorem \ref{4.3}), $\psi(N_i)=E_{h_i},$ for some  $h_1, h_2,...,
h_n\in G.$ Consider the following graph $D$ in ${\cal G}_n'(G):$\\
\centerline{\psfig{figure=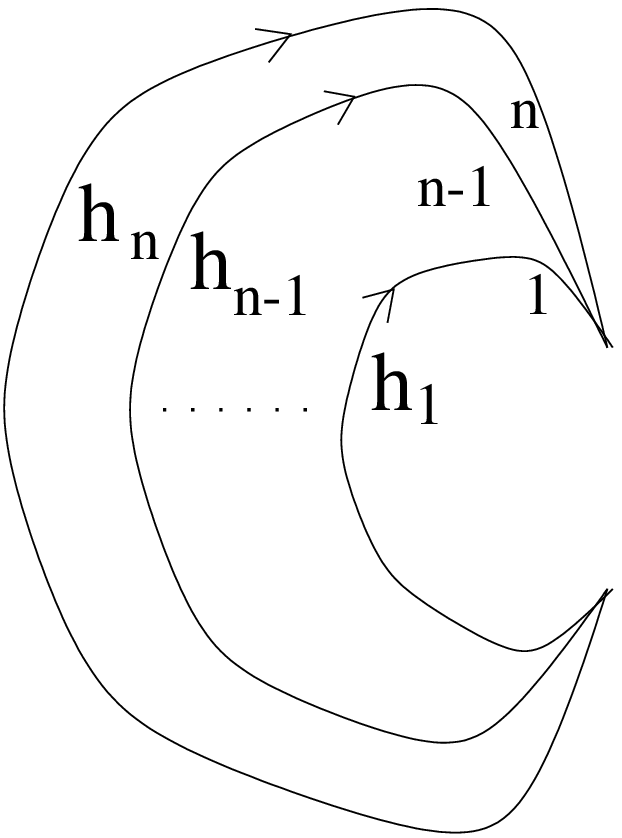,height=1.3in} \hspace*{.3in}
\psfig{figure=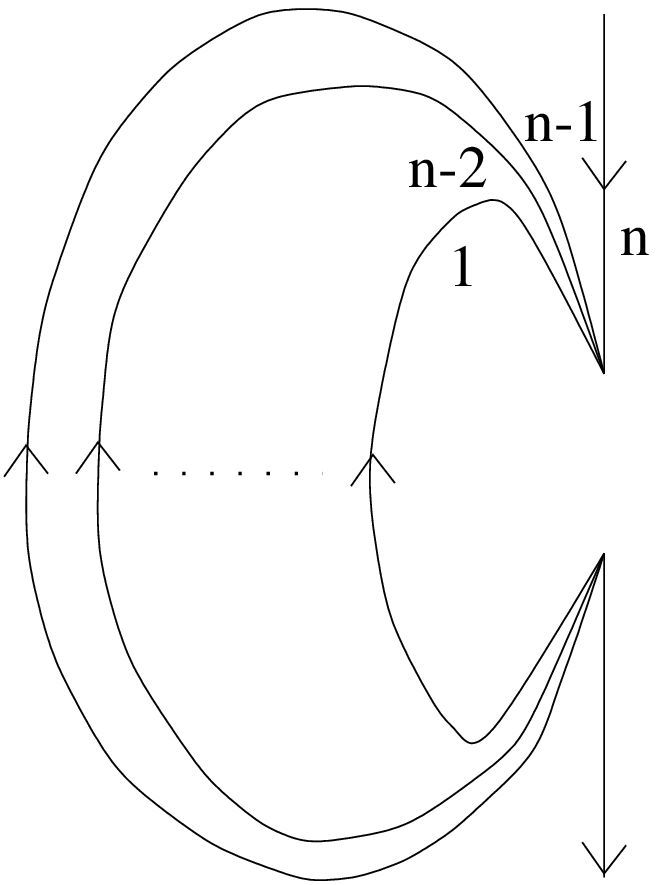,height=1.3in}}

in which we omitted labels of edges labeled by the identity in $G.$
Notice that $D$ can be also presented in the following way.

\centerline{\psfig{figure=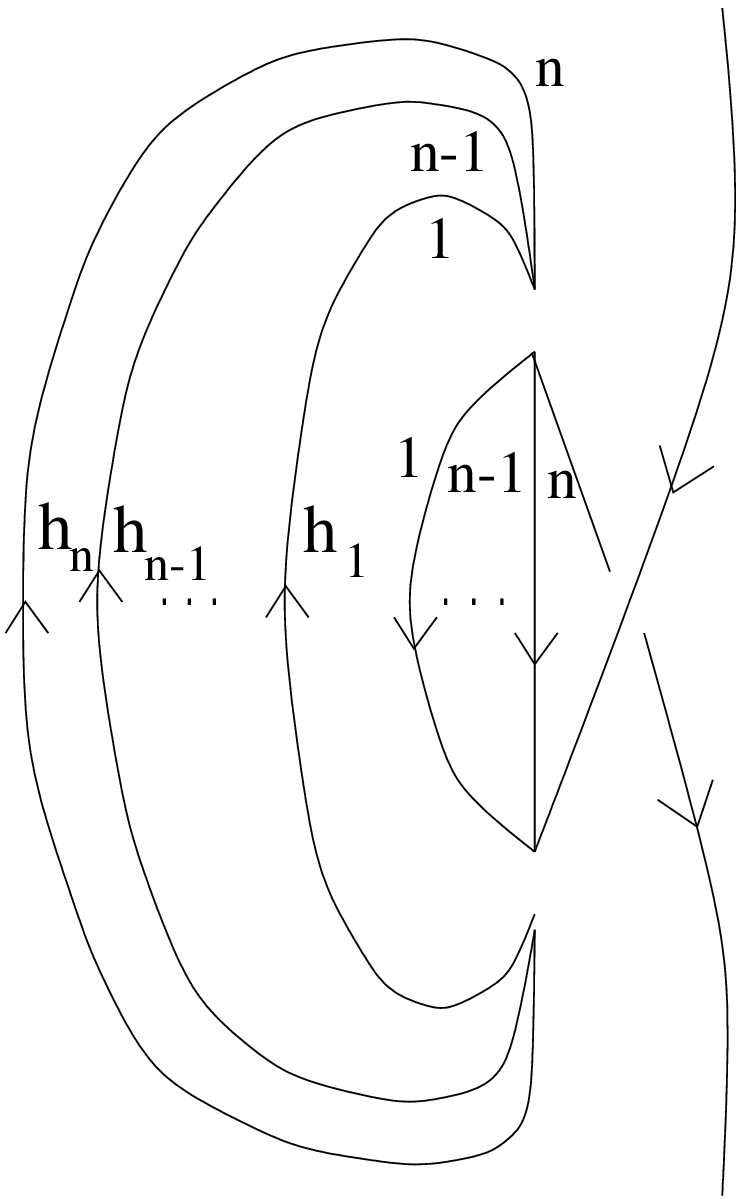,height=1.7in}}

Since the vertices of $D$ can be resolved in the two possible ways
(corresponding to the two diagrams above) we obtain the following equation:
\begin{equation}
\sum_{\sigma,\tau\in S_n} \epsilon(\sigma)\epsilon(\tau)\parbox{1.2in}
{\psfig{figure=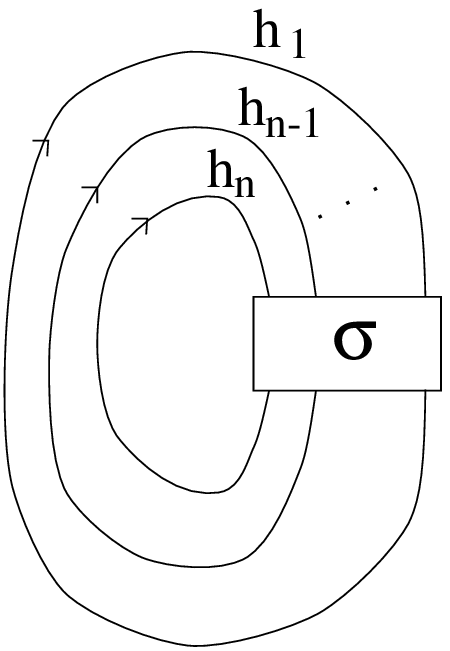,height=1.2in}}\
\parbox{1.2in}{\psfig{figure=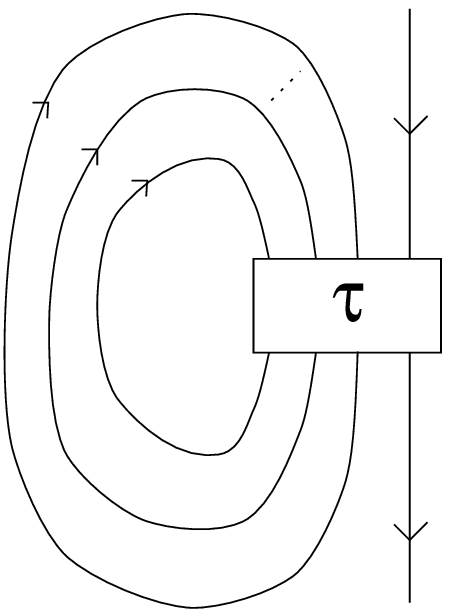,height=1.2in}}=
\sum_{\sigma,\tau\in S_n} \epsilon(\sigma)\epsilon(\tau) D_{\sigma,\tau},
\label{e16}
\end{equation}
where $D_{\sigma,\tau}$ is a graph of the form\\
\centerline{\psfig{figure=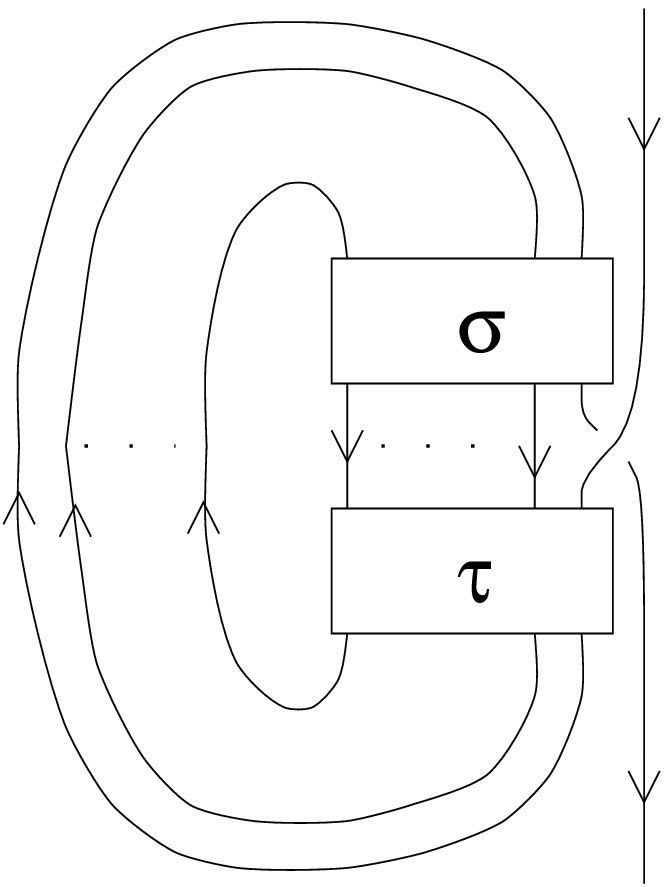,height=1.2in}}

If $\tau\in S_n$ decomposes into $k=c(\tau)$ cycles, then
$$\parbox{.8in}{\psfig{figure=diag16.eps,height=.8in}}\hspace*{.3in}=
\hspace*{.3in} \parbox{.7in}{\psfig{figure=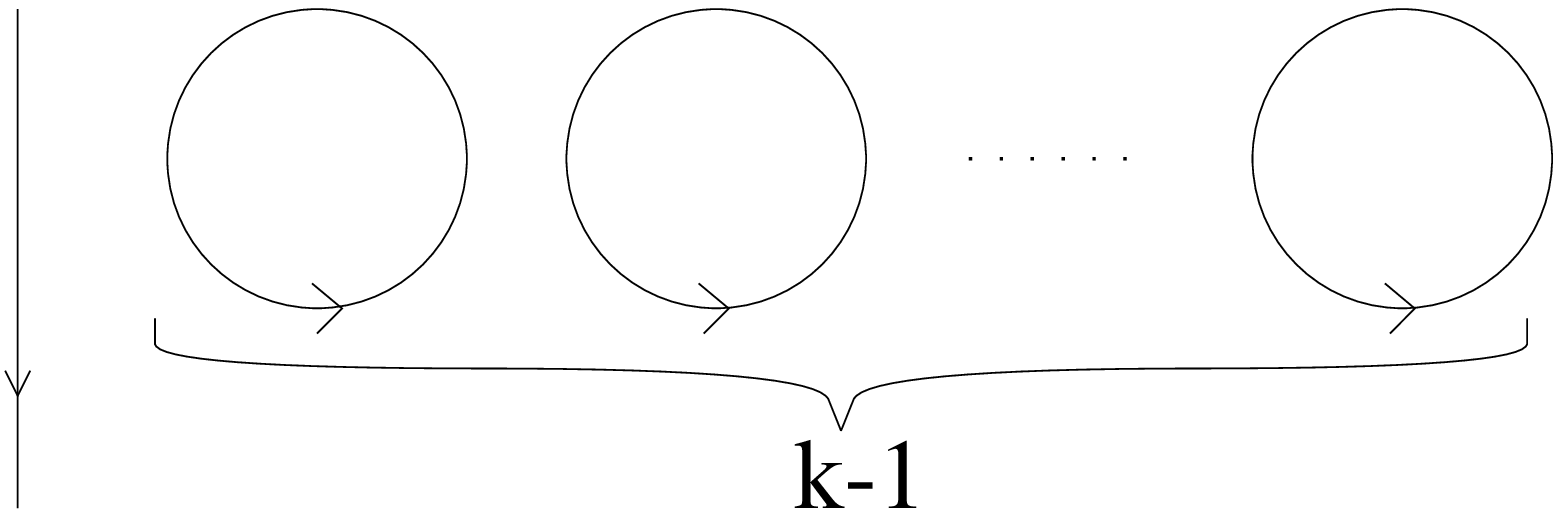,height=.7in}}$$

Therefore, by Corollary \ref{4.5},\\
$$\sum_{\tau\in S_n} \epsilon(\tau) \parbox{.8in}
{\psfig{figure=diag16.eps,height=.8in}}= \sum_{\tau\in S_n}
\epsilon(\tau) n^{c(\tau)-1}= (n-1)!$$

Notice moreover that\\
\centerline{\psfig{figure=diag9.eps,height=1.2in}}
is equal to $\psi(\Phi_{\sigma}(N_1,N_2,..,N_n)).$
Therefore the left side of (\ref{e16}) is equal to the value of $\psi$ on
$$(n-1)!\sum_{\sigma\in S_n} \epsilon(\sigma)\Phi_{\sigma}(N_1,N_2,...,N_n).$$

Now we are going to consider the right side of (\ref{e16}).
Notice that the single arc in $D_{\sigma,\tau}$ is labeled by an element
$h_{i_s}...h_{i_1}h_{i_0}\in G$ where $i_0=\tau(n),
i_1=\tau\sigma(i_0), ..., i_s=\tau\sigma(i_{s-1}),$ and
$\sigma(i_s)=n.$ Since $\tau(\sigma(i_s))=\tau(n)=i_0,$ $(i_s,i_{s-1},
..., i_1,i_0)$ is a cycle of the permutation $(\tau\sigma)^{-1}\in
S_n.$

Notice that every loop in $D_{\sigma,\tau}$ is labeled by the conjugacy
class of $h_{j_t}h_{j_{t-1}}...h_{j_1}h_{j_0},$ where
$(j_t,j_{t-1},...,j_1,j_0)$ is a cycle of
$(\tau\sigma)^{-1}\in S_n$ disjoint from $(i_s,i_{s-1},...,i_1,i_0).$
Therefore $D_{\sigma,\tau}$ is the value of $\psi$ on
$$N_{i_s}...N_{i_1}N_{i_0}Tr(N_{j_t}...N_{j_1}N_{j_0})...
Tr(N_{k_v}...N_{k_1}N_{k_0}),$$ where
$$(i_s,...,i_1,i_0)(j_t,...,j_1,j_0)...(k_v,...,k_1,k_0)$$
is the cycle decomposition of $(\tau\sigma)^{-1}.$
The above expression is equal to
$$\Phi_{(\tau\sigma)^{-1},i_s}(N_1,N_2,...,N_n)=
\Phi_{(\tau\sigma)^{-1},\sigma^{-1}(n)}(N_1,N_2,...,N_n).$$
Therefore, the right side of (\ref{e16}) is the value of $\psi$ on
$$\sum_{\sigma,\tau\in S_n} \epsilon(\sigma)\epsilon(\tau)
\Phi_{(\tau\sigma)^{-1},\ \sigma^{-1}(n)}(N_1,N_2,...,N_n)\in
T\{X_i\}_{i\in I}.$$
Let us replace $(\tau\sigma)^{-1}$ by $\kappa$ in the expression above.
Then we get $$\sum_{\sigma,\kappa\in S_n} \epsilon(\kappa)
\Phi_{\kappa,\sigma^{-1}(n)}(N_1,N_2,...,N_n)=(n-1)!
\hspace*{-.2in}\sum_{\kappa\in S_n,\ i\in \{1,2,...,n\}}
\hspace*{-.2in}\epsilon(\kappa)\Phi_{\kappa,i}(N_1,N_2,...,N_n).$$
After comparing the above algebraic descriptions of two sides of (\ref{e16})
we see that for any monomials $N_1,N_2,...,N_n$ the following element of
$T\{X_i\}_{i\in I}$ belongs to $Ker\ \psi:$

$$\sum_{\sigma\in S_n}\epsilon(\sigma)\Phi_{\sigma}(N_1,N_2,...,N_n)-
\sum_{\kappa\in S_n,\ i\in \{1,2,...,n\}}\epsilon(\kappa)
\Phi_{\kappa,i}(N_1,N_2,...,N_n).$$

By Lemma \ref{4.2} the above expression is equal to $G(N_1,N_2,...,N_n).$
Therefore $\psi(G(N_1,N_2,...,N_n))=0.$
\end{proof}

\begin{lemma}\label{4.7}\ \\
Let $N_1,N_2,...,N_{n+1}$ be any monomials in the variables $X_i,\ i\in I.$
Then\\ $\psi(F(N_1,N_2,...,N_{n+1}))=0.$
\end{lemma}

\begin{proof}
By definition, $F(N_1,N_2,...,N_{n+1})=Tr(G(N_1,N_2,...,N_n)N_{n+1}).$
By (\ref{e15}), $\psi$ commutes with the trace function.
Therefore
$$\psi(F(N_1,N_2,...,N_{n+1}))=\psi(Tr(G(N_1,N_2,...,N_n)N_{n+1}))=$$
$$Tr(\psi(G(N_1,N_2,...,N_n))\psi(N_{n+1}))=0.$$
\end{proof}

Lemmas \ref{4.6} and \ref{4.7} and the Second Fundamental Theorem
imply that the kernel of $\psi: T\{X_i\}_{i\in I}\to {\mathbb
A}_n(X,x_0)$ contains the kernel of $\pi: T\{X_i\}_{i\in I}\to
M_n(P_n(I))^{GL_n(R)}.$ Therefore we have the following corollary.

\begin{corollary}\label{4.8}\ \\
There exists an $R$-algebra homomorphism $\psi': M_n(P_n(I))^{GL_n(R)}\to
{\mathbb A}_n(X,x_0),$ such that $\psi'(A_i)=E_{g_i}$ and
$\psi'(Tr(A_{i_1}A_{i_2}...A_{i_k}))=EL_{g_{i_1}g_{i_2}...g_{i_k}},$
for any $i_1,i_2,...,i_k\in I.$
\end{corollary}

The epimorphism $\eta: P_n(I)\to Rep_n^R(G)$ introduced in Section 2 induces
an epimorphism $M_n(\eta): M_n(P_n(I))\to M_n(Rep_n^R(G))$
and, therefore, by restriction, a homomorphism $M_n(\eta)^{GL_n(R)}:
M_n(P_n(I))^{GL_n(R)} \to M_n(Rep_n^R(G))^{GL_n(R)}.$
Our goal is to show that $\psi'$ descends to
$$\psi'':M_n(Rep_n^R(G))^{GL_n(R)}\to {\mathbb A}_n(X,x_0)$$ such that the
following diagram commutes.
\begin{equation}
\begin{tabular}{ccc}
$M_n(P_n(I))^{GL_n(R)}$&&\\
$\Big \downarrow \scriptstyle{M_n(\eta)^{GL_n(R)}}$ & $\stackrel{\psi'}
{\searrow}$ &\\
$M_n(Rep_n^R(G))^{GL_n(R)}$ & $\stackrel{\psi''}{\longrightarrow}$ &
${\mathbb A}_n(X,x_0)$\\
\label{e17}
\end{tabular}
\end{equation}

In order to prove this fact we need to show that $Ker\ M_n(\eta)^{GL_n(R)}
\subset Ker\ \psi'.$ We will use the following lemma.

\begin{lemma}\label{4.9}\ \\
\begin{enumerate}
\item $Det(A_i)\in P_n(I)^{GL_n(R)}\subset M_n(P_n(I))^{GL_n(R)}$
\item $\psi'(Det(A_i))=1,$ for any $i\in I.$
\end{enumerate}
\end{lemma}

\begin{proof}
\begin{enumerate}
\item By Proposition \ref{4.4}, $Det(A_i)$ can be expressed as a linear
combination of traces of powers of $A_i.$ By Lemma \ref{1.1}(2),
$A_i^k\in M_n(P_n(I))^{GL_n(R)}$ and hence, by Proposition \ref{1.2},
$Tr(A_i^k)\in P_n(I)^{GL_n(R)}.$ Finally, by Lemma \ref{1.1}(1)
there is a natural embedding $P_n(I)^{GL_n(R)}\subset
M_n(P_n(I))^{GL_n(R)}.$

\item If $c_1, c_2, ...,c_k$ are the lengths of all cycles of
$\sigma\in S_n,$ then $\psi'$ maps\\ $Tr(A_i^{c_1})Tr(A_i^{c_2})...
Tr(A_i^{c_k})$ to a graph\\
\centerline{\psfig{figure=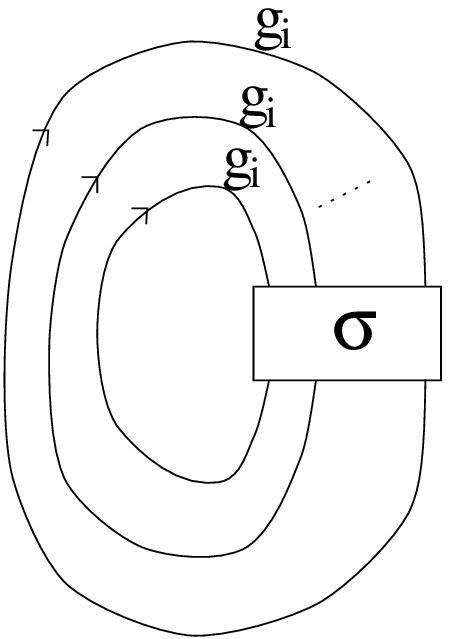,height=1.2in}}

Therefore, by Proposition \ref{4.4} and Fact \ref{3.1},
$$\psi'(Det(A_i))= {1\over n!} \sum_{\sigma\in S_n}
\epsilon(\sigma)\parbox{1in}{\psfig{figure=diage1.eps,height=1in}}={1\over n!}
\parbox{1in}{\psfig{figure=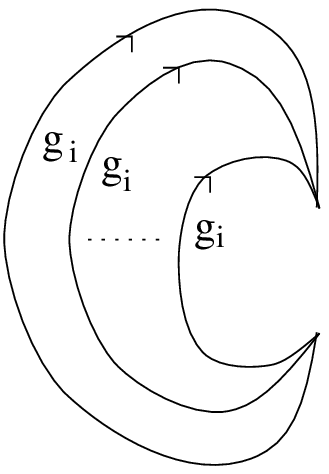,height=1in}}.$$
Analogously, $$1=\psi'(Det({\bf 1}))={1\over n!}
\parbox{1in}{\psfig{figure=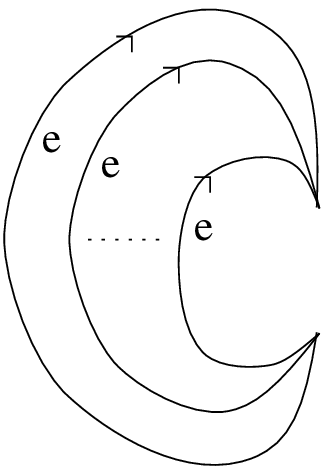,height=1in}},$$
where $e$ is the identity in $G.$
But, by (\ref{e8}) (or, equivalently, (\ref{e9}))
$$\parbox{1in}{\psfig{figure=diage2.eps,height=1in}}=
\parbox{1in}{\psfig{figure=diage3.eps,height=1in}}.$$
Therefore $\psi'(Det(A_i))=1.$
\end{enumerate}
\end{proof}

The next proposition is due to Procesi. Since the proof of this proposition
is hidden in the proof of Theorem 2.6 in \cite{Pro-2} we will recall
it here for completeness of this paper\footnote{Compare also
Proposition 9.5 \cite{B-H}.}.

\begin{proposition}\label{4.10}\ \\
Let ${\cal J}\triangleleft M_n(P_n(I))^{GL_n(R)}$
be a two sided ideal and let ${\cal J}'$ be the ideal in $P_n(I)$
generated by the entries of elements of
$M_n(P_n(I)){\cal J}M_n(P_n(I))\triangleleft M_n(P_n(I)).$
Then
\begin{enumerate}
\item $M_n(P_n(I)){\cal J}M_n(P_n(I))=M_n({\cal J}')\triangleleft M_n(P_n(I)).$
\item There is a unique $GL_n(R)$-action on $M_n(P_n(I)/{\cal J}')$ such that
the natural projection $i: M_n(P_n(I))\to M_n(P_n(I)/{\cal J}')$ is
$GL_n(R)$-equivariant.

\item $i$ induces a homomorphism
$$j:M_n(P_n(I))^{GL_n(R)}/{\cal J}\to M_n(P_n(I)/{\cal J}')^{GL_n(R)}$$
which is an isomorphism of $R$-algebras.
\end{enumerate}
\end{proposition}

\begin{proof}
\begin{enumerate}
\item This follows from the basic algebraic fact, that every
ideal ${\cal I}$ in $M_n(R),$ for any ring $R$ with $1$, is of the
form $M_n({\cal I}'),$ where ${\cal I}'$ is the ideal in $R$
generated by the entries of a generating set of the ideal ${\cal
I}.$
\item Let $B\in GL_n(R)$ and $B*$ denote the action of $B$ on
$M_n(P_n(I)).$ $B$ leaves $M_n({\cal J}')$
invariant. Indeed, any element $C\in M_n({\cal J}')$ is of the form
$\sum_i M_iC_iN_i,$ where $M_i,N_i\in M_n(P_n(I)),\ C_i\in {\cal J},$
and therefore $$B*C= \sum_i (B*M_i)(B*C_i)(B*N_i)\in M_n({\cal J}').$$
This implies that the action
of $GL_n(R)$ on $M_n(P_n(I)/{\cal J}')$ is well defined. All other
statements of (2) are obvious consequences of this fact.
\item
For any rational action of $GL_n(R)$ on any $R$-vector space $N$ there exists
a linear projection, called the Reynolds operator,
${\cal r}: N\to N^{GL_n(R)},$ with the following properties
\begin{enumerate}
\item[(i)] ${\cal r}(x)=x,$ for $x\in N^{GL_n(R)}$ and, therefore,
${\cal r}$ is an epimorphism.
\item[(ii)] ${\cal r}$ is natural with respect to $GL_n(R)$-equivariant maps
$N\to N'$
\item[(iii)] If $N$ is an algebra then ${\cal r}(xy)=x{\cal r}(y)$ and
${\cal r}(yx)={\cal r}(y)x$ for $x\in N^{GL_n(R)}$ and $y\in N.$
\end{enumerate}
For more information about this operator see \cite{MFK} or a more
elementary text \cite{Fog}.

The homomorphism $i$ restricted to $M_n(P_n(I))^{GL_n(R)}$
induces a homomorphism
$$j:M_n(P_n(I))^{GL_n(R)}/{\cal J}\to M_n(P_n(I)/{\cal
J}')^{GL_n(R)}.$$
By Property (i) of ${\cal r},$ $j$ is an epimorphism.
It remains to be proved that $j$ is injective.

Choose $i_0\in I.$ For any monomial $m$ in $P_n(I)=R[x^i_{jk}\ i\in I,
j,k=1,2, ...,n]$ we define the degree of $m$ to be the number of
appearances of the variables $x^{i_0}_{jk},$
$j,k\in\{1,2,...,n\},$ in $m.$ This induces a grading on $P_n(I).$
We can extend this grading on $M_n(P_n(I))$ as follows.
For any matrix $A=(a_{jk})\in M_n(P_n(I))$ with a single non-zero entry
$a_{st},$ $deg(A)=deg(a_{st}).$ Note that the degree of
the matrix $A_{i_0}=(x^{i_0}_{jk})\in M_n(P_n(I))$
considered in Section 2 is $1.$

Let $B\in GL_n(R).$ By the definition of the $GL_n(R)$-action on
$M_n(P_n(I))$ and by Lemma \ref{1.1}(2),
$$B\left(
\begin{matrix}
B*x_{11}^{i_0} & B*x_{12}^{i_0} & ... & B*x_{1n}^{i_0}\\
\vdots & \vdots & ... & \vdots\\
B*x_{n1}^{i_0} & B*x_{n2}^{i_0} & ... & B*x_{nn}^{i_0}
\end{matrix}
\right)B^{-1} = \left(
\begin{matrix}
x_{11}^{i_0} & x_{12}^{i_0} & ... & x_{1n}^{i_0}\\
 \vdots & \vdots & ... & \vdots\\
x_{n1}^{i_0} & x_{n2}^{i_0} & ... & x_{nn}^{i_0}
\end{matrix}\right).$$
Therefore $B*x^{i_0}_{jk}$ is
a linear combination of the variables $x^{i_0}_{j'k'},\ j',k'=1,2,...,n,$
and hence the action of $GL_n(R)$ preserves the grading of $P_n(I).$
For any $M\in M_n(P_n(I)),$ $B*M$ is a matrix obtained by applying the
action of $B$ to all entries of $M$ and then by conjugating the
resulted matrix by $B.$ Therefore the action of $GL_n(R)$ also
preserves the grading of $M_n(P_n(I)).$
The naturality of the Reynolds operators ${\cal r}:P_n(I)\to
P_n(I)^{GL_n(R)}, {\cal r}: M_n(P_n(I))\to M_n(P_n(I))^{GL_n(R)}$
implies that they also preserve the gradings.
This fact will be an important element of the proof of Proposition
\ref{4.10}(3).

We need to show that
$${\cal J}=M_n(P_n(I)){\cal J}M_n(P_n(I))\cap
M_n(P_n(I))^{GL_n(R)}.$$
However, it is sufficient to show that
$${\cal J}\supset M_n(P_n(I)){\cal J}M_n(P_n(I))\cap
M_n(P_n(I))^{GL_n(R)},$$ since the opposite inclusion is obvious.

Let $c=\sum_i a_ic_ib_i\in M_n(P_n(I))^{GL_n(R)},$ where $a_i,b_i\in
M_n(P_n(I)),c_i\in {\cal J}.$ We will show that $c\in {\cal J}.$
Since $c$ involves only finitely variables $x^i_{jk}$ and $I$ is
infinite we can choose $i_0\in I$ such that
$x^{i_0}_{jk},\ j,k=1,2,...,n,$ do not appear in $a_i,b_i,c_i.$
Thus $deg\ a_i=deg\ b_i=deg\ c_i=0.$

Consider $Tr(cA_{i_0}).$ By our assumptions about $c$ and by Lemma
\ref{1.1}(2), $cA_{i_0}\in M_n(P_n(I))^{GL_n(R)}.$ Proposition
\ref{1.2} states that $Tr: M_n(P_n(I))\to P_n(I)$ is
$GL_n(R)$-equivariant and therefore $Tr(cA_{i_0})\in M_n(P_n(I))^{GL_n(R)}.$
Thus
$$Tr(cA_{i_0})=Tr\left({\cal r}(cA_{i_0})\right)=
Tr\left(\sum_i {\cal r}(a_ic_ib_iA_{i_0})\right)=$$
$$Tr\left(\sum_i {\cal r}(b_iA_{i_0}a_ic_i)\right)=
Tr\left(\sum_i {\cal r} (b_iA_{i_0}a_i)c_i\right).$$
Note that $b_iA_{i_0}a_i$ has degree $1$ and, since ${\cal r}$
preserves the grading, ${\cal r}(b_iA_{i_0}a_i)$ is also of degree $1.$
By The First Fundamental Theorem of Invariant Theory (Theorem \ref{4.1}),
$M_n(P_n(I))^{GL_n(R)}$ is generated by the elements $A_i$ and $Tr(M),$
where $M$ varies over the set of monomials composed of non-negative
powers of matrices $A_i,
i\in I.$ By our definition of degree, $$deg(A_i)=
\begin{cases}
1& \text{if $i=i_0$}\\ 0& \text{otherwise}\\
\end{cases}$$ and
$$deg(Tr(M))={\rm \ number\ of\ appearances\ of\ } A_{i_0} {\rm\ in\ }
M.$$ Therefore, ${\cal r}(b_iA_{i_0}a_i)$ can be presented as
$$\sum_j p_{ij}A_{i_0}q_{ij}+ \sum_k Tr(s_{ik}A_{i_0})t_{ik},$$ for
some elements $p_{ij},q_{ij},s_{ik}, t_{ik}\in M_n(P_n(I))^{GL_n(R)}$
of degree $0.$
Thus $$Tr(cA_{i_0})=Tr\left(\sum_i\sum_j p_{ij}A_{i_0}q_{ij}c_i+
\sum_i \sum_k Tr(s_{ik}A_{i_0})t_{ik}c_i\right)=$$
$$Tr\left(\sum_i\sum_j p_{ij}A_{i_0}q_{ij}c_i\right)+
\sum_i \sum_k Tr(s_{ik}A_{i_0})Tr(t_{ik}c_i)=$$
$$Tr\left(\left(\sum_i\sum_j q_{ij}c_ip_{ij}+
\sum_i \sum_k Tr(t_{ik}c_i)s_{ik}\right)A_{i_0}\right).$$

Therefore $$Tr\left(\left[c-\left(\sum_i\sum_j q_{ij}c_ip_{ij}+
\sum_i \sum_k Tr(t_{ik}c_i)s_{ik}\right)\right] A_{i_0}\right)=0$$
in $M_n(P_n(I)).$ The expression in brackets above has degree $0.$
Note that if $deg\ d=0,\ d\in M_n(P_n(I)),$
then $Tr(dA_{i_0})=0$ if and only if $d=0.$ Therefore
$$c=\sum_i\sum_j q_{ij}c_ip_{ij}+ \sum_i \sum_k Tr(t_{ik}c_i)s_{ik},$$
and  hence $c\in {\cal J}.$
\end{enumerate}
This completes the proof of Proposition \ref{4.10}.
\end{proof}

Let ${\cal J}\triangleleft M_n(P_n(I))^{GL_n(R)}$ be the ideal generated by
elements $Det(A_i)-1,\ i\in I,$ and elements $A_{i_1}A_{i_2}...A_{i_k}-1$
corresponding to all defining relations $r_j= g_{i_1}g_{i_2}...g_{i_k}$ of $G.$
By Lemma \ref{4.9}(1) and Lemma \ref{1.1}(2) $Det(A_i)-1$ and
$A_{i_1}A_{i_2}...A_{i_k}-1$ are indeed elements of $M_n(P_n(I))^{GL_n(R)}$
and therefore ${\cal J}$ is well defined.
By Proposition \ref{4.10}(1) the ideal
$M_n(P_n(I)){\cal J}M_n(P_n(I))\triangleleft M_n(P_n(I))$ is equal to
$M_n({\cal J}'),$ where ${\cal J}'\triangleleft P_n(I)$ is the ideal generated
by coefficients of matrices belonging to ${\cal J}.$
Notice that ${\cal J}'$ is exactly the kernel of the epimorphism
$\eta: P_n(I)\to Rep_n^R(G)$ introduced in Section 2. Therefore by
Proposition \ref{4.10} the homomorphism
$$M_n(\eta)^{GL_n(R)}:M_n(P_n(I))^{GL_n(R)}\to M_n(Rep_n^R(G))^{GL_n(R)}$$
considered in diagram (\ref{e17}) descends to an isomorphism

$$j: M_n(P_n(I))^{GL_n(R)}/{\cal J}\to M_n(Rep_n^R(G))^{GL_n(R)}.$$

\begin{proposition}\label{4.11}\ \\
$M_n(Rep_n^R(G))^{GL_n(R)}$ is generated by the elements $j_{G,n}(g_i)$ and
$Tr\left(j_{G,n}(g_{i_1}g_{i_2}...g_{i_k})\right),$ where
$i,i_1,i_2,...,i_k\in I.$
\end{proposition}

\begin{proof}
From the paragraph preceding Proposition \ref{4.11} follows that
$M_n(\eta)^{GL_n(R)}$ is an epimorphism.
By Theorem \ref{4.1} $M_n(P_n(I))^{GL_n(R)}$ is generated by the
elements $A_i$ and $Tr(A_{i_1}A_{i_2}...A_{i_k}),$ where
$i,i_1,i_2,..., i_k\in I.$ The homomorphism $M_n(\eta)^{GL_n(R)}$
carries these elements to
$j_{G,n}(g_i)$ and $Tr\left(j_{G,n}(g_{i_1}g_{i_2}...g_{i_k})\right),$
respectively.
\end{proof}

This proposition and Theorem \ref{2.6} imply that $\Theta$ is an
epimorphism. We are going to show that it is also a monomorphism.

We have shown in Lemma \ref{4.9} that $Det(A_i)-1\in Ker\ \psi'.$
\footnote{Recall that the map $\psi'$ was defined in Corollary \ref{4.8}.}
Moreover, by the definition of $\psi',\ A_{i_1}A_{i_2}...A_{i_k}-1\in
Ker\ \psi',$
for any $i_1,i_2,...,i_k$ such that $g_{i_1}g_{i_2}...g_{i_k}=e$ in $G.$
Therefore ${\cal J}\subset Ker\ \psi'$ and we can factor $\psi'$ to
$$\psi'':M_n(Rep_n^R(G))^{GL_n(R)}\to {\mathbb A}_n(X,x_0),$$ such that
diagram (\ref{e17}) commutes and, by Corollary \ref{4.8},
\begin{itemize}
\item $\psi''(j_{G,n}(g_i))=E_{g_i}$
\item $\psi''(Tr(j_{G,n}(g_{i_1}g_{i_2}...
g_{i_k}))=EL_{g_{i_1}g_{i_2}...g_{i_k}},$ for any $i_1,i_2,...,i_k\in I.$
\end{itemize}

Recall that our assumptions about the presentation of $G$ (stated in the
paragraph following Theorem \ref{4.3}) say that the inverse
of any generator of $G$ is also a generator and that every element of
$G$ is a product of non-negative powers of generators. Thus, by Proposition
\ref{2.5}(4), ${\mathbb A}_n(X,x_0)$ is generated by the elements $E_{g_i}$ and
$EL_{g_{i_1}g_{i_2}...g_{i_k}},$ for $i,i_1,i_2,...,i_k\in I.$ Since
$\psi''\circ\Theta$ is the identity on the generators of
${\mathbb A}_n(X,x_0),$ it also is the identity
on ${\mathbb A}_n(X,x_0).$ Therefore $\Theta$ is a monomorphism.

In order to complete the proof of Theorem \ref{2.7} we need to show that
$\theta$ is also an isomorphism.

Fact \ref{2.4} implies that we have an embedding $\imath_*: {\mathbb A}_n(X)\to
{\mathbb A}_n(X,x_0),\ \imath_*(L_g)=EL_g,$ for $g\in G.$
Therefore we can consider ${\mathbb A}_n(X)$ as a subring of ${\mathbb A}_n(X,x_0).$
Moreover, by Theorem \ref{2.6}, $\theta$ is just the
restriction of $\Theta: {\mathbb A}_n(X,x_0)\to M_n(Rep_n^R(G))^{GL_n(R)}$ to
${\mathbb A}_n(X).$ Therefore $\theta$ is a monomorphism.

In order to show that $\theta$ is an epimorphism we use once again an argument
from invariant theory. By the naturality of the Reynolds operators
${\cal r}:M_n(Rep_n^R(G))\to M_n(Rep_n^R(G))^{GL_n(R)},$
${\cal r}':Rep_n^R(G)\to Rep_n^R(G)^{GL_n(R)}$
the following diagram commutes

\begin{center}
\begin{tabular}{ccc}
$M_n(Rep_n^R(G))$ & $\stackrel{Tr}{\longrightarrow}$ & $Rep_n^R(G)$\\
$\Big \downarrow \scriptstyle{{\cal r}}$ & & $\Big \downarrow \scriptstyle{
{\cal r}'}$\\
$M_n(Rep_n^R(G))^{GL_n(R)}$ &$\stackrel{Tr}{\longrightarrow}$ &
$Rep_n^R(G)^{GL_n(R)}$\\
\end{tabular}
\end{center}

Since $Tr: M_n(Rep_n^R(G))\to Rep_n^R(G)$ and all Reynolds operators
are epimorphic, $Tr: M_n(Rep_n^R(G))^{GL_n(R)}\to Rep_n^R(G)^{GL_n(R)}$ is also
epimorphic. But now commutativity of (\ref{e5}) implies that $\theta$ is
an epimorphism as well.

Therefore we have shown that $\theta$ is an isomorphism.
This completes the proof of Theorem \ref{2.7}.


\section{$SL_n$-character varieties}

In this section we present one of possible applications of Theorem \ref{2.7}
to a study of $SL_n$-character varieties.

Let $X$ be a path connected topological space whose fundamental group,
$G=\pi_1(X),$ is finitely generated. Let
$K$ be an algebraically closed field of characteristic $0.$
Recall that we noticed in Section 2 that the set of all
$SL_n(K)$-characters of $G,$ denoted by $X_n(G),$ is
an algebraic set whose coordinate ring is $Rep_n^R(G)^{GL_n(K)}/\sqrt 0.$

Let $\chi_g=Tr(j_{G,n}(g))\in Rep_n^R(G)^{GL_n(K)}/\sqrt 0,$ for any $g\in G.$
It is not difficult to check that $\chi_g$ considered as an element of
$K[X_n(G)]$ is a function which assigns to a character $\chi$ the value
$\chi(g).$ By Proposition \ref{2.5}(3) and Theorem \ref{2.7}
$Rep_n^R(G)^{GL_n(K)}$ is generated by the elements $Tr(j_{G,n}(g)).$ Therefore
the functions $\chi_g$ generate $K[X_n(G)].$

By an $SL_n$-trace identity for $G$ we mean a polynomial function
in variables $\chi_g,g\in G,$ which is identically equal to $0$ on $X_n(G).$
For example, $$\chi_g\chi_h=\chi_{gh} +\chi_{gh^{-1}},$$ is the
famous Fricke-Klein $SL_2$-trace identity valid for any group $G$ and any
$g,h\in G.$ From the above discussion it follows that the coordinate ring
of $X_n(G)$ can be considered as the quotient of the ring of polynomials in
formal variables $\chi_g, g\in G,$ by the ideal of all $SL_n$-trace identities
for $G.$ Therefore Theorem \ref{2.7} implies the following corollary.

\begin{corollary}\label{5.2}\ \\
There is an isomorphism $\Lambda: {\mathbb A}_n(X)/\sqrt 0 \to
K[X_n(G)],$ such that
$\Lambda(L_g)=\chi_g.$ Under this isomorphism the identities on graphs in $X$
induced by skein relations correspond to $SL_n$-trace identities for $G.$
\end{corollary}

The above corollary is very useful in the study of trace identities
since it makes possible to interpret them geometrically.
Consider for example the following $SL_3$-trace identity which holds for
any $\gamma_0,\gamma_1,\gamma_2,\gamma_3\in G$ and any $\chi\in
X_3(G),$ where $G$ is an arbitrary group.
\vspace*{-.1in}
$$\chi(\gamma_1)\chi(\gamma_2)\chi(\gamma_3)-\chi(\gamma_1)
\chi(\gamma_2\gamma_3)-\chi(\gamma_2)\chi(\gamma_1\gamma_3)-
\chi(\gamma_3)\chi(\gamma_1\gamma_2)+$$
$$\chi(\gamma_1\gamma_2\gamma_3)+
\chi(\gamma_1\gamma_3\gamma_2)-
\chi(\gamma_1\gamma_0)\chi(\gamma_2\gamma_0)\chi(\gamma_3\gamma_0)+$$
$$\chi(\gamma_1\gamma_0)\chi(\gamma_2\gamma_0\gamma_3\gamma_0)+
\chi(\gamma_2\gamma_0)\chi(\gamma_1\gamma_0\gamma_3\gamma_0)+
\chi(\gamma_3\gamma_0)\chi(\gamma_1\gamma_0\gamma_2\gamma_0)-$$
\begin{equation}
\chi(\gamma_1\gamma_0\gamma_2\gamma_0\gamma_3\gamma_0)-
\chi(\gamma_1\gamma_0\gamma_3\gamma_0\gamma_2\gamma_0)=0.
\label{e18}
\end{equation}

Our theory provides the following interpretation of
this identity. Let $x_0\in X$ and $G=\pi_1(X,x_0).$
Let $\gamma_0$ be a path in $X$ representing a non-trivial element of
$\pi_1(X,x_0).$
We assume that $\gamma_0$ goes along a ``hole'' in $X$ presented on
the picture below. Consider the following two, obviously
equivalent, graphs $\Gamma$ and $\Gamma'$ in $X:$
\vspace*{.2in}\\
\centerline{\psfig{figure=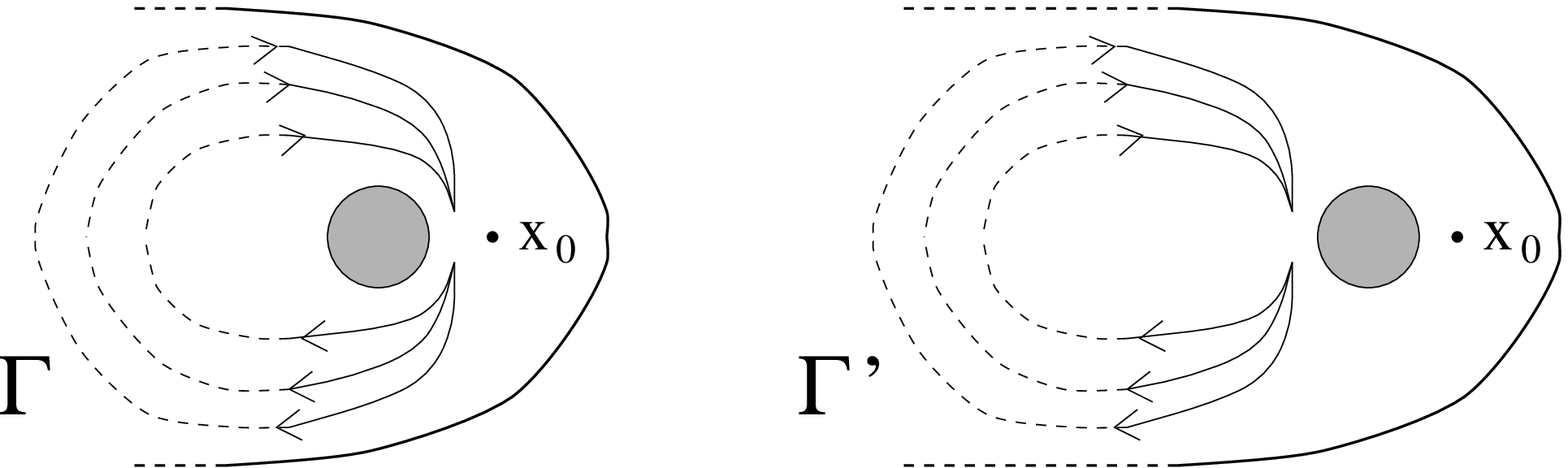,height=1in}}

The graph $\Gamma'$ is obtained from $\Gamma$ by
pulling its vertices along the ``hole'' in $X.$ The obvious resolution
of vertices in $\Gamma$ and $\Gamma'$ gives an equation involving
closed loops in $X.$ This equation corresponds
to the trace identity (\ref{e18}).

There is a large body of literature about $SL_2$-character varieties
and their applications.
However very little is known about $SL_n$-character varieties, for
$n>2.$ The reason for that is that the $SL_n$-trace
identities, like (\ref{e18}), are intractable by classical (algebraic)
methods.
Since our theory often gives a simple geometric interpretation to
complicated trace identities, it can be applied to a more detailed
study of character varieties. This idea was already used in \cite{PS-2} and
\cite{PS-3} to study $SL_n$-character varieties for $n=2.$ A
generalization of these results for $n>2,$ which is based on our skein
method, will appear in future papers. In this work we test our method
on the simplest non-trivial example -- we study $SL_3$-character
variety of the free group on two generators, $F_2=<g_1,g_2>.$ The
basic problem is to determine the minimal
dimension of the affine space in which $X_3(F_2)$ is embedded, or
equivalently, the minimal number of generators of $K[X_3(F_2)].$
A result of Procesi (Theorem 3.4(a) \cite{Pro-1}) implies that
$K[X_3(F_2)]$ is generated by the elements $\chi_{g_{i_1}g_{i_2}...g_{i_j}},$
where $j\leq 7$ and $i_1,i_2,...,i_j\in\{1,2\}.$ A direct calculation
shows that after identifying words in $g_1,g_2$ which are related by cyclic
permutations we get a set of $57$ generators of $K[X_3(F_2)].$
It is difficult to obtain any further reduction of this set in any
simple algebraic manner. However, our geometric method allows us to
reduce this problem to the study $3$-valent graphs in the twice-punctured
disc. By playing on pictures of such graphs one can reduce the number
of generators of $K[X_3(F_2)]$ to nine! These are:
$$\chi_{g_1},\ \chi_{g_2},\ \chi_{g_1^2},\ \chi_{g_2^2},\
\chi_{g_1g_2},\ \chi_{g_1^2g_2},\ \chi_{g_1g_2^2},\ \chi_{g_1^2g_2^2},
\ \chi_{g_1^2g_2^2g_1g_2}.$$
Moreover, it is possible to show that this is the minimal
number of generators and $X_3(F_2)\subset K^9$ is a solution set of
one irreducible polynomial.


\end{document}